\documentclass[11pt,a4paper,fleqn]{article}
\usepackage{amsfonts,amsmath,latexsym,amssymb,euscript,graphicx,mathrsfs}

\usepackage{graphicx}
\usepackage{color}
\usepackage{amssymb}
\usepackage{amssymb}
\usepackage[T1]{fontenc}
\usepackage{latexsym}
\usepackage{xypic}
\usepackage{eufrak}
\usepackage{euscript}
\usepackage{amsfonts,amsmath}
\usepackage{verbatim}
\usepackage[english]{babel}
\usepackage{mathrsfs}

\newtheorem{prop}{Proposition}[section]
\newtheorem{proposition}[prop]{Proposition}
\newtheorem{cor}[prop]{Corollary}
\newtheorem{lemma}[prop]{Lemma}
\newtheorem{rem}[prop]{Remark}

\newtheorem{thm}[prop]{Theorem}

\newtheorem{example}[prop]{Example}

\renewcommand{\geq}{\geqslant}
\def\leq{\leqslant}

\newcommand{\N}{\mathbb{N}}
\newcommand{\Z}{\mathbb{Z}}
\newcommand{\R}{\mathbb{R}}

\def\esp{E}

\def\e{\varepsilon}

\def\1{{\mathbf{1}}}

\def\1{{\mathbf{1}}}
\def\0.5{{\frac{1}{2}}}

\def\H{\EuFrak H}

\newcommand{\qed}{\nopagebreak\hspace*{\fill}
{\vrule width6pt height6ptdepth0pt}\par}

\begin{document}

\begin{center}
{\Large{\bf Stein's method and stochastic analysis of Rademacher functionals}}
\\~\\ by Ivan Nourdin\footnote{Laboratoire de Probabilit\'es et Mod\`eles Al\'eatoires, Universit\'e
 Pierre et Marie Curie (Paris VI), Bo\^ite courrier 188, 4 place Jussieu, 75252 Paris Cedex 05, France. Email: \texttt{ivan.nourdin@upmc.fr}},
 Giovanni Peccati\footnote{Equipe Modal'X, Universit\'{e} Paris Ouest -- Nanterre la D\'{e}fense, 200 Avenue de la République, 92000 Nanterre, and LSTA, Universit\'{e} Paris VI, France. Email: \texttt{giovanni.peccati@gmail.com}} and Gesine Reinert\footnote{Department of Statistics, University of Oxford, 1 South Parks Road, Oxford OX1 3TG, UK. Email \texttt{reinert@stats.ox.ac.uk}}\\ {\it Universit\'e Paris VI, Universit\'e Paris Ouest and Oxford University}\\~\\
\end{center}

{\small \noindent {\bf Abstract:} We compute explicit bounds in the Gaussian approximation of functionals of infinite Rade\-macher sequences. Our tools involve Stein's
method, as well as the use of appropriate discrete Malliavin operators. Although our approach does not require the classical use of exchangeable pairs, we employ a
 chaos expansion in order to construct an explicit exchangeable pair vector for
any random variable which depends on a finite set of Rademacher variables. Among several examples, which include random variables which depend on infinitely many
Rademacher variables, we provide three main applications: (i) to CLTs for multilinear forms belonging to a fixed chaos, (ii) to the Gaussian approximation of weighted
infinite $2$-runs, and (iii) to the computation of explicit bounds in CLTs for multiple integrals over sparse sets. This last application provides an alternate proof (and
several refinements) of a recent result by Blei and Janson.}
\\



\section{Introduction}
The connection between Stein's method (see e.g. \cite{Chen_Shao_sur, Reinert_sur, Stein_orig, Stein_book}) and the integration by parts formulae of stochastic analysis
has been the object of a recent and quite fruitful study. The papers \cite{NouPec_ptrf, NouPe08, NouPeRev} deal with Stein's method and Malliavin calculus (see e.g.
\cite{Nbook}) in a Gaussian setting; in \cite{NouVie} one can find extensions involving density estimates and concentration inequalities; the paper \cite{PSTU} is devoted
to explicit bounds, obtained by combining Malliavin calculus and Stein's method
 in the framework of Poisson measures. Note that all these references contain examples and applications that were
previously outside the scope of Stein's method, as they involve functionals of infinite-dimensional random fields for which the usual Stein-type techniques (such as
exchangeable pairs, zero-bias transforms or size-bias couplings), which involve picking a random index from a {\sl finite} set, seem inappropriate.

The aim of this paper is to push this line of research one step further, by combining Stein's method with an appropriate {\sl discrete} version of Malliavin calculus, in
order to study the normal approximation of the functionals of an infinite {\sl Rademacher sequence}. By this expression we simply mean a sequence $X=\{X_n : n\geq 1\}$ of
i.i.d. random variables, such that $P(X_1 =1)=P(X_1=-1) = 1/2 $. A full treatment of the discrete version of Malliavin calculus adopted in this paper can be found in
\cite{Privault} or \cite{PrivaultSchoutens}; Section \ref{SS : discreteMall} below provides a short introduction to the main objects and results.

The main features of our approach are the following:
\begin{itemize}
\item[(i)] We will be able to deal directly with the normal approximation of random variables possibly depending on the whole {\sl infinite} sequence $X$. In
    particular, our results apply to random elements not necessarily  having the form of partial sums.
\item[(ii)] We will obtain bounds expressed in terms of the
    kernels appearing in the {\sl chaotic decomposition} of a given square-integrable random variable. Note that every square-integrable functional of $X$ admits a
    chaotic decomposition into an orthogonal sum of multiple integrals.
\item[(iii)] We employ the chaotic expansion in order to construct an explicit exchangeable
    pair vector for any random variable which depends on a finite set of Rademacher variables. In particular, this exchangeable pair satisfies a linearity condition
    in conditional expectation {\sl exactly}, which allows the normal approximation results from Reinert and R\"ollin \cite{ReinertRoellin} to be applied.
\item[(iv)]
    In the particular case of random variables belonging to a {\sl fixed chaos}, that is, having the form of a (possibly infinite) series of multilinear forms of a
    fixed order, we will express our bounds in terms of norms of {\sl contraction operators} (see Section \ref{S : ZiggyStardust}). These objects also play a central
    role in the normal and Gamma approximations of functionals of Gaussian fields (see \cite{NouPec_ptrf, NouPe08, NouPeRev, NP}) and Poisson measures (see
    \cite{PSTU, PecTaq_Ber, PecTaq_SURV}). In this paper, we shall consider kernels defined on discrete sets, and the corresponding contraction norms are integrals
    with respect to appropriate tensor products of the counting measure.
\end{itemize}

The use of contraction operators in the normal approximation of Rademacher functionals seems to us a truly promising direction for further research. One striking application of these techniques is given in Section \ref{S : App}, where we deduce explicit bounds (as well as a new proof) for a combinatorial central limit theorem (CLT) recently proved by Blei and Janson in \cite{BleiJanson}, involving sequences of multiple integrals over finite ``sparse'' sets. In the original Blei and Janson's proof, the required sparseness condition emerges as an efficient way of re-expressing moment computations related to martingale CLTs. In contrast, our approach does not involve martingales or moments and (quite surprisingly) the correct sparseness condition stems naturally from the definition of contraction. This yields moreover an explicit upper bound on an adequate distance between probability measures.

Another related work is the paper by Mossel {\it et al.} \cite{MOO}, where the authors prove a general
invariance principle, allowing to compare the distribution of polynomial forms constructed from
different collections of i.i.d. random variables. In particular, in \cite{MOO} a crucial role is played by kernels ``with low influences'': we will point out in Section \ref{S :
FixedWienerBounds} that influence indices are an important component of our bounds for the normal approximation of elements of a fixed chaos. In Section \ref{SS :
questionBJ}, we also show that the combination of the findings of the present paper with those of \cite{MOO} yields an elegant answer to a question left open by Blei and
Janson in \cite{BleiJanson}.

A further important point is that our results do not hinge at all on the fact that $\mathbb{N}$ is an ordered set. This is one of the advantages both of Stein's method
and of the Malliavin-type, filtration-free stochastic analysis. It follows that our findings can be directly applied to i.i.d. Rademacher families indexed by an arbitrary
discrete set.

The reader is referred to Chatterjee \cite{Chatterjee_AOP} for another Stein-type approach to the normal approximation of functionals of finite i.i.d. vectors. Some
connections between our findings and the results proved in \cite{Chatterjee_AOP} (in the special case of quadratic functionals) are pointed out in Section \ref{SS :
quadratic}.

The paper is organized as follows. Section \ref{S : FrameW} contains results on multiple integrals, Stein's method and discrete Malliavin calculus (in particular, a new
chain rule is proved). Section \ref{S : Generalbounds} provides
 general bounds for the normal approximation of Rademacher functionals; it gives first examples and the exchangeable pair construction. In Section \ref{S :
 FixedWienerBounds} we are interested in
random variables belonging to a fixed chaos. Section \ref{S : 2runs} focuses on sums of single and double integrals and related examples. Section \ref{S : App} is devoted
to further applications and refinements, mainly involving integrals over sparse sets. Technical proofs are deferred to
the Appendix. 

\section{Framework and main tools}\label{S : FrameW}
\subsection{The setup}\label{setup}
As before, we shall denote by $X=\left\{ X_{n}:n\geq 1\right\} $ a standard \textsl{Rademacher sequence}. This means that $X$ is a collection of i.i.d. random variables,
defined on a probability space $(\Omega, {\mathcal{F}}, P)$ and such that $P\left[X_{1}=1\right] =P\left[ X_{1}=-1\right] =1/$2. To simplify the subsequent discussion,
for the rest of the paper we shall suppose that $\Omega = \{-1,1\}^{\mathbb{N}}$ and $P = [\frac{1}{2}\{\delta_{-1} + \delta_1\}]^{\mathbb{N}}$. Also, $X$ will be defined
as the canonical process, that is: for every $\omega = (\omega_1,\omega_2,...)\in\Omega$, $X_n (\omega)=\omega_n$. Starting from $X$, for every $N\geq 1$ one can build a
random signed measure $\mu _{(X,N)}$ on $\{1,...,N\}$, defined for every $A\subset \{1,...,N\}$ as follows:
$$
\mu _{(X,N)}\left( A\right) =\sum_{j\in A}X_{j}.
$$
It is clear that $\sigma\{X\} = \sigma\{\mu_{(X,N)}: N\geq 1\}$. We shall sometimes omit the index $N$ and write $\mu_{(X,N)} = \mu_X$, whenever the domain of integration is clear from the context.

\smallskip

We define the set $D=\left\{ \left( i_{1},i_{2}\right) \in \mathbb{N}^{2}:i_{1}=i_{2}\right\} $ to be the \textsl{diagonal} of $\mathbb{N}^{2}$. For $n\geq 2$, we put
\begin{equation}\label{simplex1}
\Delta _{n}=\left\{ \left( i_{1},...,i_{n}\right) \in \mathbb{N}^{n}:\text{ \ the }i_{j}\text{'s are all different}\right\},
\end{equation}
and, for $N,n\geq 2$,
\begin{equation}\label{simplex2}
\Delta _{n}^{N}=\Delta _{n}\cap \left\{ 1,...,N\right\} ^{n}
\end{equation}
(so that $\Delta _{n}^{N}=\varnothing $ if $N<n$). The random signed measures $\mu _{(X,N)}$ satisfy that, for every $A,B \subset \{1,...,N\}$,
\begin{eqnarray}
\mu _{(X,N)}^{\otimes 2}\left( \left[ A\times B\right] \cap D\right) =\sum_j X_j^{2}\mathbf{1}_{\{j\in A\}} \mathbf{1}_{\{j \in B\}}\label{control} =\kappa \left( A \cap
B\right) =\sharp \left\{ j:j\in A\cap B\right\},
\end{eqnarray}
where $\kappa $ is the counting measure on $\mathbb{N}$. The application $A\mapsto \mu _{(X,N)}^{\otimes 2}\left( \left[ A\times A\right] \cap D\right) $ is called the
\textsl{ diagonal measure } associated with $\mu_{(X,N)}$. The above discussion simply says that, for every $N\geq 1$, the diagonal measure associated with $\mu_{(X,N)}$
is the restriction on $\{1,...,N\}$ of the counting measure $\kappa$. Note that, for finite sets $A,B$ one has also $$E[\mu_X (A) \mu_X (B)] = \kappa(A\cap B).$$
\begin{rem}
{\rm In the terminology e.g. of \cite{PecTaq_SURV} and \cite{Sur}, the collection of random variables $$ \{\mu_{(X,N)} (A) : A \subset \{1,...,N\}, \, N \geq 1 \} $$
defines an \textsl{independently scattered random measure} (also called a \textsl{completely random measure}) on $\mathbb{N}$,
 with \textsl{control measure} equal to $\kappa$. In particular, if $A_1,...,A_k$  is a collection of mutually disjoint finite subsets of $\mathbb{N}$, then the random
 variables $\mu_X (A_1),...,\, \mu_X (A_k)$ are mutually independent. }
\end{rem}

\subsection{The star notation}\label{S : ZiggyStardust}
For $n\geq 1$, we denote by $\ell ^{2}\left( \mathbb{N}\right) ^{n}$ the class of the kernels (= functions) on $\mathbb{N}^{n}$
that are square integrable with respect to $\kappa ^{\otimes n}$; $\ell ^{2}\left( \mathbb{N} \right)^{\circ n}$
is the subset of $\ell ^{2}\left( \mathbb{N}\right) ^{n}$ composed of symmetric kernels; $\ell _{0}^{2}\left( \mathbb{N}\right) ^{n}$
is the subset of $\ell ^{2}\left( \mathbb{N}\right) ^{n}$ composed of kernels vanishing on diagonals, that is, vanishing on the complement
of $\Delta _{n}$ ; $\ell _{0}^{2}\left( \mathbb{N}\right) ^{\circ n}$ is the subset of $\ell _{0}^{2}\left( \mathbb{N}\right) ^{n}$ composed of symmetric kernels.

\smallskip

For every $n,m\geq 1$, every $r=0,...,n\wedge m$, every $l=0,...,r$, and every $f\in \ell _{0}^{2}\left( \mathbb{N}\right)^{\circ n}$ \ and \ $g\in \ell _{0}^{2}\left(
\mathbb{N}\right) ^{\circ m}$, we denote by $f\star _{r}^{l}g$ the function of $n+m-r-l$ variables obtained as follows: $r$ variables are identified and, among these, $l$
are integrated out with respect to the counting measure $\kappa $. Explicitly,
\begin{eqnarray*}
&&f\star _{r}^{l}g\left( i_{1},...,i_{n+m-r-l}\right) \\ &=&\sum_{a_{1},...,a_{l}}f\left( i_{1},...,i_{n-r};i_{n-r+1},...,i_{n-l};a_{1},...,a_{l}\right) \\
 &&\text{ \ \ \ \ \ \ \ \ \ \ \ \ \ \ \ \ \ \ \ \ \ \ \ \ \ \ \ \ }\times g\left( i_{n-l+1},...,i_{n+m-r-l};i_{n-r+1},...,i_{n-l};a_{1},...,a_{l}\right) ;
\end{eqnarray*}
note that, since $f$ and $g$ vanish on diagonals by assumption, the sum actually runs over all vectors $\left( a_{1},...,a_{l}\right) \in \Delta _{l}.$ For instance,%
\begin{equation*}
f\star _{0}^{0}g\left( i_{1},...,i_{n+m}\right) =f\otimes g\left( i_{1},...,i_{n+m}\right) =f\left( i_{1},...,i_{n}\right) g\left( i_{n+1},...,i_{n+m}\right) ,
\end{equation*}
and, for $r\in\{1,\ldots,n\wedge m\}$,
\begin{eqnarray*}
&&f\star _{r}^{r}g\left( i_{1},...,i_{n+m-2r}\right) \\ &=&\sum_{a_{1},...,a_{r}}f\left( i_{1},...,i_{n-r};a_{1},...,a_{r}\right) g\left(
i_{n-r+1},...,i_{n+m-2r};a_{1},...,a_{r}\right) .
\end{eqnarray*}
In particular, if $r=m=n$
\begin{equation*}
f\star _{m}^{m}g=\sum_{a_{1},...,a_{m}}f\left( a_{1},...,a_{m}\right) g\left( a_{1},...,a_{m}\right) = \langle f,g \rangle_{\ell ^{2}\left( \mathbb{N}\right) ^{\otimes
m}} .
\end{equation*}

\begin{example}{\rm
If $n=m=2$, one has
\begin{eqnarray*}
 &&f\star _{0}^{0}g \left( i_{1},i_2, i_3,i_{4} \right)
=f \left( i_{1},i_2 \right) g\left( i_3, i_4\right) \,\, ; \,\, f\star _{1}^{0}g\left( i_{1},i_2, i_3 \right) =f\left( i_{1},i_2\right) g\left( i_3, i_2\right) \, ;\\ &&
f\star _{1}^{1}g\left( i_{1},i_2 \right) =\Sigma_a  f\left( i_{1}, a\right) g\left( i_2, a\right) \,\, ; \,\, f\star _{2}^{0}g\left( i_{1},i_2 \right) = f\left(
i_{1},i_2\right) g\left( i_1, i_2\right) \, ;\\ && f\star _{2}^{1}g\left( i \right) = \Sigma_a f\left( i, a \right) g\left( i, a \right) \,\, ; \,\, f\star _{2}^{2}g =
\Sigma_{a_1, a_2} f\left( a_1, a_2 \right) g\left( a_1, a_2\right).
\end{eqnarray*}
}
\end{example}

\smallskip
\begin{rem}
{\rm In general, a kernel of the type $f\star _{r}^{l}g$ may be not symmetric and may not vanish on diagonal sets. }
\end{rem}

The next two lemmas collect some useful facts for subsequent sections. Their proofs are postponed
in the appendix.
\begin{lemma} \label{L : estimates}
Let $f\in \ell_0 ^{2}\left( \mathbb{N}\right) ^{\circ n}$, $g\in \ell_0 ^{2}\left( \mathbb{N}\right) ^{\circ m}$ $(n,m\geq 1)$ and $0\leq l\leq r\leq n\wedge m$. Then:
\begin{enumerate}
\item $f\star_r^lg\in \ell ^{2}\left( \mathbb{N}\right) ^{\otimes (n+m-r-l)}$ and $\| f\star_r^l g\|_{\ell ^{2}\left( \mathbb{N}\right) ^{\otimes (n+m-r-l)}}
 \leq
\| f\|_{\ell ^{2}\left( \mathbb{N}\right) ^{\otimes n}} \| g\|_{\ell ^{2}\left( \mathbb{N}\right) ^{\otimes m}}; $
\item if $n\geq 2$ then
\begin{eqnarray*}
\max_{j\in\mathbb{N}}  \left[\sum_{(b_1,...,b_{n-1}) \in\mathbb{N}^{n-1}} \!\!\!\!\!\!\! f^2(j,b_1,...,b_{n-1})\right]^2  & \leq & \|f\star _n^{n-1}
f\|^2_{\ell^2(\mathbb{N})}\\ &\leq & \|f\|^2_{\ell^2(\mathbb{N})^{\otimes n}}\!\times\!  \max_{j\in\mathbb{N}} \sum_{(b_1,...,b_{n-1}) \in\mathbb{N}^{n-1}}
\!\!\!\!\!\!\!f^2(j,b_1,...,b_{n-1})
\end{eqnarray*}
and, if $l=1,\ldots,n$,
\begin{eqnarray}
\|f\star _l^{l-1} g\|^2_{\ell^2(\mathbb{N})^{\otimes (n+m -2l +1)}} &=& \sum_{j=1}^{\infty} \| f(j,\cdot)\star_{l-1}^{l-1} g(j,\cdot)\|_{\ell ^{2}\left(
\mathbb{N}\right) ^{\otimes (n+m-2l)}} ^2 \notag\\ &\leq& \|f\star _n^{n-1} f\|_{\ell^2(\mathbb{N})} \times \|g\|^2_{\ell^2(\mathbb{N})^{\otimes m}}.\label{MODALX}
\end{eqnarray}
\item $ \| f\star_1^0 f \|_{\ell^2(\mathbb{N})^{\otimes (2n-1)}} = \| f\star_n^{n-1} f \|_{\ell^2(\mathbb{N})} $ and, for every $l=2,...,n$ $(n\geq 2)$,
$$
\| f\star_l^{l-1} f \|_{\ell^2(\mathbb{N})^{\otimes (2n-2l+1)}} \leq \|f\star_{l-1}^{l-1} f \times \mathbf{1}_{\Delta^c_{2(n-l+1)}} \|_{\ell^2(\mathbb{N}) ^{\otimes
(2n-2l+2)} }\leq \|f\star_{l-1}^{l-1} f\|_{\ell^2(\N)^{\otimes (2n-2l+2)}}.
$$
Here, $\Delta^c_q$ stands for the complement of the set $\Delta^c_q$, that is, $\Delta^c_q$ is the collection of all vectors $(i_1,...,i_q)$ such that $i_k=i_l$ for
at least one pair $(l,k)$ such that $l\neq k$.
\end{enumerate}
\end{lemma}
\begin{rem}\label{after L : estimates}
{\rm
\begin{itemize}
\item[1.] According e.g. to \cite{MOO}, the quantity $$\sum_{(b_1,...,b_{n-1})\in\N^{n-1}} \!\!\!\!\! f^2(j,b_1,...,b_{n-1})$$ is called the {\sl influence of the
    $j^{\rm th}$ coordinate} on $f$. \item[2.] When specializing Lemma \ref{L : estimates} to the case $n=2$, one gets
\begin{eqnarray}
\max_{j\in\N} \left[\sum_{i\in\N} f^2(i,j)\right]^2 & \leq & \|f\star_2^1 f\|^2_{\ell^2 (\mathbb{N})}\label{traceequality}\\ &=&\|f\star_1^1 f\times{\bf
1}_{\Delta^c_2}\|^2_{\ell^2 (\mathbb{N})^{\otimes 2}} \leq\|f\star_1^1 f\|^2_{\ell^2 (\mathbb{N})^{\otimes 2}}  ={\bf Trace}\{[f]^4\}.\notag
\end{eqnarray}
Here, $[f]$ denote the infinite array $\{f(i,j):\,i,j\in\N\}$ and $[f]^4$ stands for the fourth power of $[f]$ (regarded as the kernel of a Hilbert-Schmidt
operator).
\end{itemize}
}
\end{rem}

\begin{lemma}\label{L : convLemma}
Fix $n,m\geq 1$, and let $f_k\in \ell^2_0(\mathbb{N})^{\circ n}$ and $g_k\in \ell^2_0(\mathbb{N})^{\circ m}$, $k\geq 1$, be such that $f_k
\underset{k\to\infty}{\longrightarrow} f$ and $g_k \underset{k\to\infty}{\longrightarrow} g$, respectively, in $\ell^2_0(\mathbb{N})^{\circ n}$ and
$\ell^2_0(\mathbb{N})^{\circ m}$. Then, for every $r=0,...,n\wedge m$, $f_k \star_r^r g_k \underset{k\to\infty}{\longrightarrow} f \star_r^r g$ in
$\ell^2(\mathbb{N})^{\otimes m+n-2r}$.
\end{lemma}

\subsection{Multiple integrals, chaos and product formulae}

Following e.g. \cite{PrivaultSchoutens}, for every $q\geq 1$ and every $f\in \ell _{0}^{2}\left( \mathbb{N}\right) ^{\circ q}$ we denote by $J_{q}\left( f\right) $ the \textsl{multiple integral} (of order $q$) of $f$ with respect to $X$. Explicitly,
\begin{eqnarray}\label{MWII}
J_q(f) &=& \sum_{(i_1,...,i_q)\in\mathbb{N}^q} f(i_1,...,i_q) X_{i_1}\cdot\cdot\cdot X_{i_q}=\sum_{(i_1,...,i_q)\in \Delta^q} f(i_1,...,i_q) X_{i_1}\cdot\cdot\cdot
X_{i_q}\\ &=&q!\sum_{i_1<\ldots<i_q} f(i_1,...,i_q) X_{i_1}\cdot\cdot\cdot X_{i_q},\notag
\end{eqnarray}
where the possibly infinite sum converges in $L^2(\Omega)$. Note that, if $f$ has support contained in some finite set $\{1,...,N\}^q$, then $J_{q}\left(f\right) $
is a true integral with respect to the product signed measure $\mu _{(X,N)}^{\otimes q}$, that is,
\begin{equation*}
J_{q}\left( f\right) =\int_{\{1,...,N\}^q}f \, d\mu _{(X,N)}^{\otimes q}=\int_{\Delta^N _{q}}f \, d\mu _{(X,N)}^{\otimes q}.
\end{equation*}
One customarily sets $\ell ^{2}\left( \mathbb{N}\right)^{\circ 0} = \mathbb{R}$, and $J_0(c) = c$, $\forall c\in \mathbb{R}$. It is known (see again
\cite{PrivaultSchoutens}) that, if $f \in \ell ^{2}_0\left( \mathbb{N} \right)^{\circ q}$ and $g \in \ell ^{2}_0\left( \mathbb{N}
 \right)^{\circ p}$, then one has the isometric relation
\begin{eqnarray} \label{iso}
E[ J_q(f) J_p(g)] = {\bf 1}_{\{q=p\}} q! \langle f, g \rangle_{ {\ell^2(\mathbb{N})^{ \otimes q} } } .
\end{eqnarray}
In particular,
\begin{eqnarray} \label{iso2}
E [J_q(f)^2]  =  q! \|f\|^2_{\ell^2(\mathbb{N})^{ \otimes q}}.
\end{eqnarray}
The collection of all random variables of the type $J_n(f)$, where $f \in \ell ^{2}_0\left( \mathbb{N}\right)^{\circ q} $, is called the $q$th \textsl{chaos} associated
with $X$.

\begin{rem} {\rm Popular names for random variables of the type $J_q(f)$ are \textsl{Walsh chaos} (see e.g. \cite[Ch. IV]{LedTal}) and \textsl{Rademacher chaos} (see e.g.
\cite{BleiJanson, delapena, KW}). In Meyer's monograph \cite{Mey92} the collection $\{J_q(f) : f\in
  \ell ^{2}_0\left( \mathbb{N}\right)^{\circ q}
 , \, q\geq 0 \} $
is called the \textsl{toy Fock space} associated with $X$. }
\end{rem}

Recall (see \cite{PrivaultSchoutens}) that one has the following {\sl chaotic decomposition}: for every $F\in L^2(\sigma\{X\})$ (that is, for every square integrable
functional of the sequence $X$) there exists a unique sequence of kernels $f_n \in
  \ell ^{2}_0\left( \mathbb{N}\right)^{\circ n}
 $, $n\geq 1$, such that
\begin{equation}\label{Chaos}
F=E(F) + \sum_{n\geq 1} J_n (f_n) = E(F) + \sum_{n\geq 1} n! \sum_{ i_1 <i_2<...<i_n} f_n(i_1,...,i_n)X_{i_1}\cdot\cdot\cdot X_{i_n},
\end{equation}
where the second equality follows from the definition of $J_n(f_n)$, and the series converge in $L^2$.

\begin{rem}
{\rm Relation (\ref{Chaos}) is equivalent to the statement that the set
\begin{equation}\label{basis}
\{1\}\cup \bigcup_{n\geq 1}\{ X_{i_1}\ldots X_{i_n}:\,1\leq i_1<\ldots<i_n\}
\end{equation}
is an orthonormal basis of $L^2(\sigma\{X\})$. An immediate proof of this fact can be deduced from basic harmonic analysis. Indeed, it is well-known that the set $\Omega
= \{-1,+1\}^\mathbb{N}$, when endowed with the product structure of coordinate multiplication, is a compact Abelian group (known as the \textsl{Cantor group}), whose
unique normalized Haar measure is the law of $X$. Relation (\ref{Chaos}) then follows from the fact that the dual of $\Omega$ consists exactly in the mappings
$\omega\mapsto 1$ and $\omega\mapsto X_{i_1}(\omega)\cdot\cdot\cdot X_{i_n}(\omega)$, $i_n>...>i_1\geq 1$, $n\geq 1$. See e.g. \cite[Section VII.2]{BleiBook}. }
\end{rem}

Given a kernel $f$ on $\mathbb{N}^n$, we denote by $\widetilde{f}$ its canonical symmetrization, that is:

\begin{equation}\notag
\widetilde{f}\left( i_{1},...,i_{n}\right) =\frac{1}{n!}\sum_{\sigma} f (i_{\sigma ( 1 )},...,i_{\sigma (n)} ),
\end{equation}
where $\sigma $ runs over the $n!$ permutations of $\left\{ 1,...,n\right\} $. The following result is a standard multiplication formula between multiple integrals of
possibly different orders. Since we did not find this result in the literature (see, however, \cite[formula (15)]{PrivaultSchoutens}), we provide a simple combinatorial
proof in Section \ref{A : ProofProduct}.

\begin{proposition}\label{Product formula}
For every $n,m\geq 1$, every $f\in \ell _{0}^{2}\left( \mathbb{N}\right) ^{\circ n}$ and $g\in \ell _{0}^{2}\left( \mathbb{N}\right)^{\circ m}$, one has that%
\begin{eqnarray}\label{TheProductFormula}
&& J_{n}\left( f\right) J_{m}\left( g\right) = \sum_{r=0}^{n\wedge m}r!\binom{n}{r}\binom{m}{r}J_{n+m-2r}\left[ \left( \widetilde{f\star _{r}^{r}g}\right)
\mathbf{1}_{\Delta _{n+m-2r}}\right] . \label{s}
\end{eqnarray}
\end{proposition}

\begin{rem}{\rm
Proposition \ref{Product formula} yields that the product of two multiple integrals is a linear combination of square-integrable random variables. By induction, this
implies in particular that, for every $f\in \ell _{0}^{2}\left( \mathbb{N}\right) ^{\circ n}$ and every $k\geq 1$, $E|J_n(f)|^k<\infty$, that is, the elements belonging
to a given chaos have finite moments of all orders. This fact can be alternatively deduced from standard hypercontractive inequalities, see e.g. Theorem 3.2.1 in
\cite{delapena} or \cite[Ch. 6]{KW}. Note that similar results hold for random variables belonging to a fixed Wiener chaos of a Gaussian field (see e.g. \cite[Ch.
V]{Janson}). }
\end{rem}
\subsection{Finding a chaotic decomposition}\label{SS : how?}
The second equality in (\ref{Chaos}) clearly implies that, for every $i_1<...<i_n$,
$$
n!f_n(i_1,...,i_n) = E(F \times X_{i_1}\cdot\cdot\cdot X_{i_n}).
$$
In what follows, we will point out two alternative procedures yielding an explicit expression for the kernels $f_n$.

\medskip

\noindent(i) {\it M\"{o}bius inversion (Hoeffding decompositions)}. We start with the following observation: if $F = F(X_1,...,X_d)$ is a random variable that depends
uniquely on the first $d$ coordinates of the Rademacher sequence $X$, then necessarily $F$ admits a finite chaotic decomposition of the form
\begin{equation}\label{finitechaos}
F=E(F) +\sum_{n=1}^d \sum_{1\leq i_1<\ldots<i_n\leq d } n! f_n (i_1 , \ldots, i_n) X_{i_1}\cdot \cdot \cdot X_{i_n}.
\end{equation}
Note that the sum in the chaotic decomposition (\ref{finitechaos}) must  stop at $d$: indeed, the terms of order greater than $d$ are equal to zero since, if not, they
would involve e.g. products of the type $X_{j_1}\cdot\cdot\cdot X_{j_{d+1}}$ with all the $j_a$'s different, which would not be consistent with the fact that $F$ is
$\sigma\{X_1,\ldots,X_d\}$-measurable. Now note that one can rewrite (\ref{finitechaos}) as $F=E(F) + \sum_{I\subset\{1,\ldots,d\}}G(I)$, where, for
$I:=\{i_1,\ldots,i_a\}$,
\begin{equation*}
G(I) = a! f_{a} (i_1 , \ldots, i_a) X_{i_1}\ldots X_{i_a}.
\end{equation*}
By exploiting the independence, one obtains also that, for every $J := \{j_1,...,j_n\}\subset \{1,\ldots,d\}$,
\begin{equation*}
H(J): = E[F-E(F) | X_{j_1},\ldots,X_{j_n}] = \sum_{I\subset J} G(I),
\end{equation*}
thus yielding that, by inversion (see e.g. \cite[p. 116]{Stan}), for every $n=1,...,d$,
\begin{equation}\label{Hoeffding}
n! f_{n} (j_1 , ..., j_n) X_{j_1}\cdot \cdot \cdot X_{j_n} = \sum_{\{i_1,\ldots,i_a\} \subset \{j_1,\ldots,j_n\}} (-1)^{n-a} E[F-E(F) | X_{i_1},...,X_{i_a}].
\end{equation}
Formula (\ref{Hoeffding}) simply means that, for a fixed $n$, the sum of all random variables of the type $n! f_{n} (j_1 , ..., j_n) X_{j_1}\cdot \cdot \cdot X_{j_n} $
coincides with the $n^{\rm th}$ term in the \textsl{Hoeffding-ANOVA} decomposition of $F$ (see e.g. \cite{Krinott}). By a density argument (which is left to the reader),
representation (\ref{Hoeffding}) extends indeed to random variables depending on the {\sl whole} sequence $X$.

\medskip

\noindent (ii) {\it Indicators expansions}. Assume that we can write $F=F(X_1,\ldots,X_d)$ with $F:\{-1,+1\}^{d}\to\R$. Taking into account all the possibilities, we
have
$$
F=\sum_{(\e_1,\ldots,\e_d)\in\{-1,+1\}^d} F(\e_1,\ldots,\e_d)\prod_{i=1}^d {\bf 1}_{\{X_i=\e_i\}}.
$$
Now, the crucial point is the identity ${\bf 1}_{\{X_i=\e_i\}}=\frac12(1+\e_iX_i)$. Hence
\begin{eqnarray}\label{chaos2}
F=2^{-d}\sum_{(\e_1,\ldots,\e_d)\in\{- 1, +1\}^d} F(\e_1,\ldots,\e_d)\prod_{i=1}^d (1+\e_iX_i).
\end{eqnarray}
But
\begin{eqnarray*}
\prod_{i=1}^d (1+\e_iX_i) &= &1 + \sum_{1\leq i_1\leq d} \e_{i_1}X_{i_1} +\sum_{1\leq i_1<i_2\leq d} \e_{i_1}\e_{i_2}X_{i_1}X_{i_2}\\ &&+\sum_{1\leq i_1<i_2<i_3\leq d}
\e_{i_1}\e_{i_2}\e_{i_3}X_{i_1}X_{i_2}X_{i_3} + \ldots +\e_{i_1}\ldots\e_{i_d}X_{i_1}\ldots X_{i_d};
\end{eqnarray*}
inserting this in \eqref{chaos2} one can deduce the chaotic expansion of $F$.

\subsection{Discrete Malliavin calculus and a new chain rule}\label{SS : discreteMall}
We will now define a set of discrete operators which are the analogues of the classic Gaussian-based Malliavin operators (see e.g. \cite{Janson, Nbook}). The reader is
referred to \cite{Privault} and \cite{PrivaultSchoutens} for any unexplained notion and/or result.

\smallskip

The operator $D$, called the \textsl{gradient operator}, transforms random variables into random sequences. Its domain, noted ${\rm dom}D$, is given by the class of
random variables $F\in L^2(\sigma\{X\})$ such that the kernels $f_n\in \ell^2_0(\N)^{\circ n}$ in the chaotic expansion $F=E(F) + \sum_{n\geq 1} J_n (f_n)$ (see
(\ref{Chaos})) verify the relation
$$
\sum_{n\geq 1} nn!\|f_n\|^2_{\ell^2(\mathbb{N})^{\otimes n}} <\infty.
$$
In particular, if $F=F(X_1,\ldots,X_d)$ depends uniquely on the first $d$ coordinates of $X$, then $F\in{\rm dom}D$. More precisely, $D$ is an operator with values in
$L^2(\Omega\times\mathbb{N}, P\otimes\kappa)$, such that, for every $F=E(F)+\sum_{n\geq 1} J_n(f_n)\in{\rm dom}D$,
\begin{eqnarray} \label{graddef}
D_k F  = \sum_{n\geq 1} n J_{n-1} (f_n(\cdot, k) ), \quad k \geq 1,
\end{eqnarray}
where the symbol $f_n(\cdot, k)$ indicates that the integration is performed with respect to $n-1$ variables. According e.g. to \cite{Privault,PrivaultSchoutens}, the
gradient operator admits the following representation. Let $\omega = (\omega_1, \omega_2, \ldots) \in \Omega$, and set
$$\omega_+^k = (\omega_1, \omega_2, \ldots, \omega_{k-1}, +1, \omega_{k+1} , \ldots )$$
and
$$\omega_{-}^k = (\omega_1, \omega_2, \ldots, \omega_{k-1}, -1, \omega_{k+1} , \ldots )$$
to be the sequences obtained by replacing the $k^{\rm th}$ coordinate of $\omega$, respectively, with $+1$ and $-1$. Write $F_k^\pm$ instead of $F(\omega_\pm^k)$ for simplicity.
Then, for every $F \in {\rm dom}D$,
 \begin{eqnarray}\label{D_as_differential}
 D_kF(\omega) = \frac{1}{2}
 \big(F_k^+-F_k^-\big)
 , \,\, k\geq 1.
 \end{eqnarray}
\smallskip
\begin{rem}\label{GWsoon}{\rm It is easily seen that, if the random variable $F\in L^2(\sigma\{X\})$ is such that the mapping $(\omega , k) \mapsto \frac{1}{2}
(F_k^+-F_k^-)(\omega)$ is an element of $L^2(\Omega\times\mathbb{N}, P\otimes\kappa)$, then necessarily $F\in{\rm dom}D$.
}
\end{rem}

\medskip

We write $\delta$ for the adjoint of $D$, also called the \textsl{divergence operator}. The domain of $\delta$ is denoted by ${\rm dom}\delta$, and is such that ${\rm
dom}\delta \subset L^2(\Omega\times\mathbb{N}, P\otimes\kappa)$. Recall that $\delta$ is defined via the following \textsl{integration by parts formula}: for every
$F\in{\rm dom}D$ and every $u\in {\rm dom}\delta$
\begin{equation}\label{Intbyparts}
E[F\delta(u)] = E[\langle DF, u \rangle _{\ell^2(\mathbb{N})}] = \langle DF , u \rangle_{L^2(\Omega\times\mathbb{N}, P\otimes\kappa)}.
\end{equation}

\smallskip

Now set $L^2_0(\sigma\{X\})$ to be the subspace of $L^2(\sigma\{X\})$ composed of centered random variables. We write $L : L^2(\sigma\{X\}) \to L^2_0(\sigma\{X\})$ for
the \textsl{Ornstein-Uhlenbeck operator}, which is defined as follows. The domain ${\rm dom}L$ of $L$ is composed of
 random variables $F=E(F) + \sum_{n\geq 1} J_n (f_n)\in L^2(\sigma\{X\})$ such that
$$
\sum_{n\geq 1} n^2n!\|f_n\|^2_{\ell^2(\mathbb{N})^{\otimes n}} <\infty,
$$
and, for $F\in{\rm dom}L$,
\begin{equation}\label{Ldef}
LF = -\sum_{n\geq 1} n J_n(f_n).
\end{equation}
With respect to \cite{Privault,PrivaultSchoutens}, note that we choose to add a minus in the right-hand side of (\ref{Ldef}), in order to facilitate the connection with
the paper \cite{NouPec_ptrf}.
One crucial relation between the operators $\delta$, $D$ and $L$ is that
\begin{equation}\label{deltaelle}
\delta D = -L.
\end{equation}

\smallskip

The inverse of $L$, noted $L^{-1}$, is defined on $F\in L^2_0(\sigma\{X\})$, and is given by
\begin{equation}\label{LInvdef}
L^{-1}F = -\sum_{n\geq 1} \frac1n J_n(f_n).
\end{equation}

\begin{lemma}\label{fund-ipp}
Let $F\in {\rm dom}D$ be centered, and $f:\R\to\R$ be such that $f(F)\in {\rm dom}D$. Then $ E\big[Ff(F)\big]=E\big[ \langle Df(F),-DL^{-1}F\rangle_{\ell^2(\N)}\big]. $
\end{lemma}
{\it Proof}. Using (\ref{deltaelle}) and (\ref{Intbyparts}) consecutively, we can write
\begin{eqnarray*}
E\big[Ff(F)\big]=E\big[LL^{-1}Ff(F)\big] =-E\big[\delta DL^{-1}Ff(F)\big] =E\big[ \langle Df(F),-DL^{-1}F\rangle_{\ell^2(\N)}\big].
\end{eqnarray*}
\qed

\smallskip

Finally, we define $\{P_t : t \geq 0\} = \{e^{tL} : t \geq 0\} $ to be the the \textsl{semi-group} associated with $L$, that is,
\begin{equation}\label{ou}P_t F = \sum_{n=0}^\infty e^{-nt} J_n(f_n), \quad t \geq 0, \quad
\mbox{for }F = E(F) + \sum_{n=1}^\infty J_n(f_n)\in L^2(\sigma\{X\}).
\end{equation}

\smallskip

The next result will be useful throughout the paper.

\begin{lemma}\label{lm1}
Let $F\in {\rm dom} D$ and fix $k\in\mathbb{N}$. Then:
\begin{enumerate}
\item The random variables $D_kF$, $D_k L^{-1}F$, $F_k^+$ and $F_k^-$ are independent of $ X_k$. \item It holds that $| F_k^+ - F|\leq 2 |D_k F|$ and $| F_k^- -
    F|\leq 2 |D_k F|$, $P$-almost surely. \item If $F$ has zero mean, then $E\|DL^{-1}F\|^2_{\ell^2(\mathbb{N})} \leq E\|DF\|^2_{\ell^2(\mathbb{N})}$ with equality
    if and only if $F$ is an element of the first chaos.
\end{enumerate}
\end{lemma}
{\it Proof}. 1. One only needs to combine the definition of $F_k^\pm$ with (\ref{D_as_differential}).\\ 2. Use $F_k^\pm-F=\pm(F_k^+-F_k^-){\bf 1}_{\{X_k=\mp 1\}}=\pm 2
D_kF {\bf 1}_{\{X_k=\mp 1\}}$.\\ 3. Let us consider the chaotic expansion of $F$: $$F=\sum_{n\geq 1} J_n(f_n).$$ Then $-D_kL^{-1}F =  \sum_{n\geq 1}
J_{n-1}\big(f_n(\cdot,k)\big)$  and $D_k F = \sum_{n\geq 1} n J_{n-1}\big(f_n(\cdot,k)\big)$. Therefore, using
the isometric relation \eqref{iso},
\begin{eqnarray*}
E\|DL^{-1}F\|^2_{\ell^2(\mathbb{N})}&=&E\sum_{k\in\N} \left( \sum_{n\geq 1} J_{n-1}\big(f_n(\cdot,k)\big) \right)^2\\ &=&\sum_{n\geq 1} (n-1)!
\|f_n\|^2_{\ell^2(\N)^{\otimes n}}\\ &\leq&\sum_{n\geq 1} n^2(n-1)! \|f_n\|^2_{\ell^2(\N)^{\otimes n}}\\ &=&E\sum_{k\in\N} \left( \sum_{n\geq 1}
nJ_{n-1}\big(f_n(\cdot,k)\big) \right)^2=E\|DF\|^2_{\ell^2(\mathbb{N})}.
\end{eqnarray*}
Moreover, the previous equality shows that we have equality if and only if $f_n=0$ for all $n\geq 2$, that is, if and only if $F$ is an element of the first chaos. \qed

We conclude this section by proving a chain rule involving deterministic functions of random variables in the domain of $D$. It should be compared with the classic chain
rule of the Gaussian-based Malliavin calculus (see e.g. \cite[Prop. 1.2.2]{Nbook}).

\begin{prop}\label{chainrule} {\rm (Chain Rule)}. Let $F\in{\rm dom}D$ and $f:\R\to\R$ be thrice differentiable with bounded third derivative. Assume moreover that $f(F)\in {\rm dom}D$. Then,
for any integer $k$, $P$-a.s.:
$$
\left| D_k f(F) - f'(F)D_k F +\frac12\big( f''(F_k^+)+f''(F_k^-) \big) (D_k F)^2 X_k \right|\leq\frac{10}{3}|f'''|_\infty |D_k F|^3.
$$
\end{prop}
{\it Proof}. By a standard Taylor expansion,
\begin{eqnarray*}
D_kf(F)&=&\frac12\big(f(F_k^+)-f(F_k^-)\big) =\frac12\big(f(F_k^+)-f(F)\big) - \frac12\big(f(F_k^-)-f(F)\big)\\ &=&\frac12 f'(F)(F_k^+-F) +\frac14 f''(F) (F_k^+-F)^2 +
R_1 \\ &&-\frac12 f'(F)(F_k^--F) -\frac14 f''(F) (F_k^--F)^2 + R_2 \\ &=&f'(F)D_k F + \\ && \frac18 \big(f''(F_k^+)+f''(F_k^-)\big) \big((F_k^+-F)^2 -(F_k^--F)^2\big) +
R_1+R_2+R_3\\ &=&f'(F)D_k F - \frac12 \big(f''(F_k^+)+f''(F_k^-)\big) (D_kF)^2X_k + R_1+R_2+R_3,
\end{eqnarray*}
where, using Lemma \ref{lm1},
\begin{eqnarray*}
|R_1|&\leq& \frac1{12}|f'''|_\infty \big|F_k^+-F\big|^3 \leq \frac23 |f'''|_\infty |D_k F|^3\\ |R_2|&\leq& \frac1{12}|f'''|_\infty \big|F_k^--F\big|^3 \leq \frac23
|f'''|_\infty |D_k F|^3\\ |R_3|&\leq& \frac1{8}|f'''|_\infty \left(\big|F_k^+-F\big|^3 +\big|F_k^--F\big|^3\right)\leq 2 |f'''|_\infty |D_k F|^3.
\end{eqnarray*}
By putting these three inequalities together, the desired conclusion follows. \qed
\subsection{Stein's method for normal approximation}
Stein's method is a collection of probabilistic techniques, using differential operators in order to assess quantities of the type
$$
\big|E[h(F)] - E[h(Z)]\big|,
$$
where $Z$ and $F$ are generic random variables, and the function $h$ is such that the two expectations are well defined. In the specific case where
$Z\sim\mathscr{N}(0,1)$, with $\mathscr{N}(0,1)$ a standard Gaussian law, one is led to consider the so-called \textsl{Stein equation} associated with $h$, which is
classically given by
\begin{equation}\label{SteinGaussEq}
h(x)-E[h(Z)] = f'(x)-xf(x), \quad x\in\R.
\end{equation}
A solution to (\ref{SteinGaussEq}) is a function $f$, depending on $h$, which is Lebesgue a.e.-differentiable, and such that there exists a version of $f'$ verifying
(\ref{SteinGaussEq}) for every $x\in\R$. The following result collects findings by Stein \cite{Stein_orig, Stein_book}, Barbour \cite{Barbour1990}, Daly \cite{Daly} and
G\"{o}tze \cite{Goetze-91}. Precisely, the proof of Point (i) (which is known as \textsl{Stein's lemma}) involves a standard use of the Fubini theorem (see e.g.
\cite[Lemma 2.1]{Chen_Shao_sur}). Point (ii) is proved e.g. in \cite[Lemma II.3]{Stein_book} (for the estimates on the first and second derivative), in \cite[Theorem
1.1]{Daly} (for the estimate on the third derivative) and in \cite{Barbour1990} and \cite{Goetze-91}
(for the alternative estimate on the first derivative). From here onwards, denote by ${\mathcal C}_b^k$ the set of all real-valued bounded functions with
bounded derivatives up to $k^{\rm th}$  order.
\begin{lemma}\label{Stein_Lemma_Gauss}
\begin{enumerate}
\item[\rm (i)] Let $W$ be a random variable. Then, $W\stackrel{\rm Law}{=}Z\sim \mathscr{N}(0,1)$ if and only if
for every continuous, piecewise continuously differentiable function $f$ such that $\mathbb{E}|f'(Z)|$ $<$ $\infty$,
\begin{equation}\label{SteinTruc}
\mathbb{E}[f'(W)-Wf(W)]=0.
\end{equation}
\item[\rm (ii)]
If $h \in {\mathcal C}_b^2$,
 then (\ref{SteinGaussEq}) has a solution $f_h$ which is thrice differentiable and such that
$\|f_h'\|_\infty \leq 4\|h\|_\infty$,
 $\|f_h''\|_\infty \leq 2 \|h'\|_\infty$ and $\|f_h'''\|_\infty \leq 2 \|h''\|_\infty$.
We also have $\| f_h' \| \leq \| h'' \|_\infty$.
\end{enumerate}
\end{lemma}

Now fix a function $h\in {\mathcal C}_b^2$,
consider a Gaussian random variable $Z\sim \mathscr{N}(0,1)$, and let
$F$ be a generic square integrable random variable. Integrating both sides of (\ref{SteinGaussEq})
with respect to the law of $F$ gives
\begin{equation}
\big|E[h(F)]-E[h(Z)]\big| = \big|E[f_h'(F)-Ff_h(F)]\big| \label{BoundH}
\end{equation}
with $f_h$ the function given in Lemma \ref{Stein_Lemma_Gauss}. In the following sections we will show that, if $F$ is a functional of the infinite Rademacher sequence
$X$, then the quantity in the RHS of (\ref{BoundH}) can be successfully assessed by means of discrete Malliavin operators.

\begin{rem}\label{R : weak}{\rm Plainly, if the sequence $F_n$, $n\geq 1$, is such that $\big|E[h(F_n)]-E[h(Z)]\big|\rightarrow 0$ for every $h\in\mathcal{C}^2_b$, then $F_n \stackrel{\rm Law}{\rightarrow} Z$. }
\end{rem}

\section{General Bounds for Rademacher functionals}\label{S : Generalbounds}
\subsection{Main bound}
The following result combines  Proposition \ref{chainrule} with Lemma \ref{Stein_Lemma_Gauss}-(ii) and Relation (\ref{BoundH}) in order to estimate expressions of the
type $|E[h(F)] - E[h(Z)]|$, where $F$ is a square integrable functional of the infinite Rademacher sequence $X$, and $Z\sim\mathscr{N}(0,1)$.

\begin{thm}\label{firstprop}
Let $F\in{\rm dom} D$ be centered and such that $\sum_k E\big|D_k F\big|^4<\infty$. Consider a function $h \in
{\mathcal C}_b^2$ and let
$Z\sim\mathscr{N}(0,1)$. Then,
\begin{eqnarray}\label{MainBound1}
&&|E[h(F)] - E[h(Z)]|\leq \min( 4 \|h\|_\infty  , \| h '' \|_\infty) B_1  + \|h''\|_\infty   B_2 ,
\end{eqnarray}
where
\begin{eqnarray}
B_1 &= & E\big|1-\langle DF,-DL^{-1}F \rangle_{\ell^2(\mathbb{N})}\big|
\leq  \label{MainBound2.5} \sqrt{E\big[(1-\langle DF, -DL^{-1}F\rangle_{\ell^2(\mathbb{N})})^2\big]} ;
\notag\\ B_2 &=& \frac{20}{3}  E\left[\left\langle |DL^{-1}F|, |DF|^3 \right\rangle_{\ell^2(\mathbb{N})}\right]. \label{MainBound3}
\end{eqnarray}
\end{thm}
{\it Proof}. Since $h\in {\mathcal C}_b^2$,
Equality (\ref{BoundH}) holds. Observe that, since the first two derivatives of
$f_h$ are bounded, one has that $f_h(F)\in{\rm dom}D$ (this fact can be proved by using a Taylor expansion, as well as the assumptions on $DF$ and the content of Remark
\ref{GWsoon}). Using Lemma \ref{fund-ipp}, one deduces that
$$
E\big[f_h'(F)-Ff_h(F)\big] = E\big[f_h'(F)-\langle Df_h(F), -DL^{-1}F\rangle_{\ell^2(\mathbb{N})}\big].
$$
We now use again the fact that $f''_h$ is bounded: as an application of the first point of Lemma \ref{lm1}, we deduce therefore that
$$
E\left[ D_kL^{-1}F\times \big(f_h''(F_k^+)+f_h''(F_k^-) \big) (D_k F)^2 X_k \right]=0, \quad k\geq 1;
$$
in particular, the boundedness of $f_h''$, along with the fact that $E|D_k F|^4<\infty$ and Lemma \ref{lm1}-(3), ensure that the previous expectation is well-defined.
Finally, the desired conclusion follows from the chain rule proved in Proposition \ref{chainrule}, as well as the bounds on $\|f'_h\|_\infty$ and $\|f'''_h\|_\infty$
stated in Lemma \ref{Stein_Lemma_Gauss}-(ii). \qed
\begin{rem} {\rm
\begin{itemize}
\item[1.] The requirement that $\sum_k E\big|D_k F\big|^4<\infty$ is automatically fulfilled whenever $F$ belongs to a {\sl finite} sum of chaoses. This can be
    deduced from the hypercontractive inequalities stated e.g. in \cite[Ch. 6]{KW}. \item[2.] Since we are considering random variables possibly depending
 on the whole sequence $X$ and having an infinite chaotic expansion, the expectation
in (\ref{MainBound3}) may actually be infinite. \item[3.] Fix $q\geq 1$ and let $F$ have the form of a multiple integral of the type $F=J_q(f)$, where
$f\in\ell_0^2(\mathbb{N})^{\circ q}$. Then, one has that
\begin{eqnarray}
\langle DF,-DL^{-1}F\rangle_{\ell^2(\mathbb{N})} = \frac 1q \|DF \|^2_{\ell^2(\mathbb{N})}, \label{boundINT1}\\ \langle |DL^{-1}F|, |DF|^3
\rangle_{\ell^2(\mathbb{N})} = \frac1q  \|DF \|^4_{\ell^4(\mathbb{N})} . \label{boundINT2}
\end{eqnarray}
\item[4.] Let $G$ be an isonormal Gaussian process (see \cite{Janson} or \cite{Nbook}) over some separable Hilbert space $\H$, and assume that $F\in L^2(\sigma\{G\})$
    is centered and differentiable in the sense of Malliavin. Then the Malliavin derivative of $F$, noted $DF$, is a $\H$-valued random element, and in
    \cite{NouPec_ptrf} the following upper bound is established (via Stein's method) between the laws of $F$ and $Z\sim \mathscr{N}(0,1)$:
    $$
    d_{TV} (F, Z) \leq
    2
    E|1- \langle DF, -DL^{-1}F\rangle_\H |,
    $$
    where $L^{-1}$ is here the inverse of the Ornstein-Uhlenbeck generator associated with $G$ and
    $d_{TV}$ is the total variation distance between the laws of $F$ and $Z$.
\item[5.] Let $N$ be a compensated Poisson measure over some measurable space
 $(A,\mathcal{A})$, with $\sigma$-finite control measure given by $\mu$.
 Assume that $F\in L^2(\sigma\{N\})$ is centered and differentiable in the sense of Malliavin.
 Then the derivative of $F$ is a random processes which a.s. belongs to $L^2(\mu)$.
 In \cite[Section 3]{PSTU} the following upper bound for the Wasserstein distance (see Section \ref{newsectionbyGesine} )
between the laws of $F$ and $Z\sim \mathscr{N}(0,1)$ is proved (again by means of Stein's method):
    $$
    d_W (F, Z) \leq E|1- \langle DF, -DL^{-1}F\rangle_{L^2(\mu)} | + E\int_A |D_a F|^2 \times |D_aL^{-1}F| \mu(da).
    $$
\end{itemize}
}
\end{rem}
\subsection{First examples: Rademacher averages}
A (possibly infinite) \textsl{Rademacher average} is just an element of the first chaos associated with $X$, that is, a random variable of the type
\begin{equation}\label{infBA}
F = \sum_{i=1}^\infty \alpha_i X_i, \quad \mbox{with $\alpha\in\ell^2(\N)$}.
\end{equation}
See e.g. \cite[Ch. IV]{LedTal} for general facts concerning random variables of the type (\ref{infBA}). The next result is a consequence of Theorem \ref{firstprop} and of
the relations (\ref{boundINT1})--(\ref{boundINT2}). It yields a simple and explicit upper bound for the Gaussian approximation of Rademacher averages.

\begin{cor}\label{C : Averages} Let $F$ have the form (\ref{infBA}), let $Z\sim\mathscr{N}(0,1)$ and consider
 $h  \in
{\mathcal C}_b^2$. Then
\begin{equation}\label{averagebound}
\big|E[h(F)]-E[h(Z)]\big|\leq \min(4\|h\|_\infty,\|h''\|_\infty) \left|1-\sum_{i=1}^\infty \alpha_i^2 \right| + \frac{20}3\|h''\|_\infty \sum_{i=1}^\infty \alpha_i^4.
\end{equation}
\end{cor}

The proof of Corollary \ref{C : Averages} (whose details are left to the reader) is easily deduced from the fact that, for $F$ as in (\ref{infBA}), one has that
(\ref{boundINT1})-(\ref{boundINT2}) with $q=1$ hold, and $D_i F =\alpha_i$ for all $i\geq 1$. We now describe two applications of Corollary \ref{C : Averages}.
The first one provides a faster than Berry-Ess\'{e}en rate for partial sums of Rademacher random variables, using smooth test functions.

\begin{example}\label{ex3.4} {\rm Let $F_n$, $n\geq 1$, be given by
$$
F_n = \frac1{\sqrt{n}}\sum_{i=1}^n X_i.
$$
Then, Relation (\ref{averagebound}) in the special case $\alpha_i= {\bf 1}_{\{i\leq n\}}\times n^{-1/2}$ yields the bound
$$
\big|E[h(F)]-E[h(Z)]\big|\leq \frac{20}{3n}\|h''\|_\infty,
$$
where $Z\sim\mathscr{N}(0,1)$, implying a faster rate than in the classical Berry-Ess\'{e}en estimates. This faster rate arises from our use of smooth test functions; a
related result is obtained in \cite{Gold_Rein_97} using a coupling approach. }\end{example}

The next example demonstrates that our techniques allow to easily tackle some problems which, due to their non-finite description, have previously escaped the Stein's
method treatment.

\begin{example}{\rm
For every $r\geq 2$, we set
$$F_r = \sqrt{r} \sum_{i \geq r} \frac{X_i}{i}. $$
The random variable $F_r$ has the form (\ref{infBA}), with $\alpha_i = {\bf 1}_{\{i\geq r\}} \sqrt{r}/i$. We have that
$$
\left\vert 1 - \sum_{i\in\N} \alpha_i^2 \right\vert = \left\vert 1 - r\sum_{i \geq r} \frac{1}{i^2} \right\vert =\left|r\sum_{i\geq
r}\left(\frac1i-\frac1{i+1}-\frac1{i^2}\right)\right| =r \sum_{i \geq r} \frac{1}{i^2(i+1) } \leq \sum_{i\geq r} \frac{1}{i(i+1)}=\frac{1}{r}
$$
and also
\begin{eqnarray*}
\sum_{i\in\N} \alpha_i^4 = \sum_{i \geq r} \frac{r^2}{i^4} \leq  \sum_{i \geq r} \frac{1}{i(i-1)} =  \frac{1}{r-1}.
\end{eqnarray*}
It follows from Relation (\ref{averagebound}) that, for $Z\sim\mathscr{N}(0,1)$,
\begin{eqnarray*}
\left| E[h(F_r)]-E[h(Z)] \right| \leq \frac{\min(4\|h\|_\infty,\|h''\|_\infty)}{r}+ \frac{20\|h''\|_\infty}{3(r-1)}.
\end{eqnarray*}
In particular, $F_r \stackrel{\rm Law}{\rightarrow} Z$ as $r \rightarrow \infty$. }
\end{example}

\begin{rem}\label{Mehler}{\rm (A Mehler-type representation) Take an independent copy of $X$, noted $X^* = \{X^*_1,X^*_2,...\}$, fix $t>0$, and
consider the sequence $X^{t} = \{X_1^{t},X^t_2,...\}$ defined as follows: $X^t$ is a sequence of i.i.d. random variables such that, for every $k\geq 1$, $X_k^t = X_k$
with probability $e^{-t}$ and $ X_k^t = X^*_k $ with probability $1-e^{-t}$ (we stress that the choice
 between $X$ and $X^*$ is made separately for every $k$, so that one can have for instance $X^t_1 = X_1, \,X^t_2 = X^*_2,...$ and so on). Then
 $X^t$ is a Rademacher sequence and one has the following representation: for every $F\in L^2(\sigma(X))$ and every $\omega = (\omega_1,\omega_2,...)\in\Omega =
 \{-1,+1\}^\mathbb{N}$
\begin{equation}\label{pitti}
P_t F(\omega) = E[F(X^t)|X = \omega],
\end{equation}
where $P_t$ is the Ornstein-Uhlenbeck semigroup given in \eqref{ou}. To prove such a representation  of $P_t$, just consider a
 random variable of the type $X_{j_1}\times\cdot\cdot\cdot\times X_{j_d}$, and then use a density argument. Now observe that,
for a centered $F=\sum_{n\geq 1}J_n(f_n)\in L^2_0(\sigma\{X\})$, Equation \eqref{LInvdef} holds,
and consequently
\begin{equation*}
-D_k L^{-1} F  =\sum_{n\geq 1} J_{n-1}(f_n(k,\cdot)) = \int_0^\infty e^{-t}P_t D_kF dt = \int_0^\infty e^{-t}E[D_kF(X^t)|X] dt,
\end{equation*}
so that
\begin{equation*}
\langle DF , -DL^{-1}F\rangle_{\ell^2(\mathbb{N})} = \int_0^\infty e^{-t}\sum_{k\geq 1} D_kF E[D_kF(X^t)|X]\, dt = \int_0^\infty e^{-t}\langle DF , E[DF(X^t)|X]
\rangle_{\ell^2(\mathbb{N})}  \, dt.
\end{equation*}
Note that the representation (\ref{pitti}) of the Ornstein-Uhlebeck semigroup is not specific of Ra\-de\-ma\-cher sequences, and could be e.g. extended to the case where
$X$ is i.i.d. standard Gaussian (see for instance \cite[Section 2.2]{MOO}). This construction would be different from the one leading to the usual \textsl{Mehler formula}
(see e.g. \cite[Formula (1.67)]{Nbook}). See Chatterjee \cite[Lemma 1.1]{Chatterjee} and Nourdin and Peccati \cite[Remark 3.6]{NouPec_ptrf} for further connections between Mehler formulae and
Stein's method in a Gaussian setting.
 }
\end{rem}

\subsection{A general exchangeable pair construction}

Remark \ref{Mehler} uses a particular example of exchangeable pair of infinite Rademacher sequences. We shall now show that the chaos decomposition approach links  in
very well with the method of \textsl{exchangeable pairs}. Let us assume that $F = F(X_1,...,X_d)$ is a random variable that depends uniquely on the first $d$ coordinates
${\mathbf{X}}_d =(X_1, \ldots, X_d)$ of the Rademacher sequence $X$, with  finite chaotic decomposition of the form \eqref{finitechaos}. Now assume that $E(F)=0$ and
$E(F^2)=1$ so that
\begin{eqnarray} \label{fdec}
F & = &  \sum_{n=1}^d \sum_{1\leq i_1<...<i_n\leq d }   n! f_n (i_1 , ..., i_n) X_{i_1}\cdot \cdot \cdot X_{i_n} =  \sum_{n=1}^d J_n(f_n).
\end{eqnarray}
A natural exchangeable pair construction is as follows. Pick an index $I$ at random, so that $P(I=i) = \frac{1}{d}$ for $i=1, \ldots, d$, independently of $X_1,...,X_d$,
and if $I=i$ replace $X_i$ by an independent copy $X_i^*$ in all sums in the decomposition \eqref{fdec}
which involve $X_i$. Call the resulting expression $F'$. Also
denote the vector of Rademacher variables with the exchanged component by ${\mathbf{X}}'_d$. Then $(F,F')$ forms an exchangeable pair.

 In \cite{ReinertRoellin} an
embedding approach is introduced, which suggests enhancing the pair  $(F,F')$ to a pair of vectors $(  \mathbf{W} ,  \mathbf{W} ')$ which then ideally satisfy the
linearity condition
\begin{eqnarray} \label{generalexch}
E ( \mathbf{W'} - \mathbf{W} | \mathbf{W} ) = - \Lambda \mathbf{W}
\end{eqnarray}
with $\Lambda$ being a deterministic matrix. In our example, choosing  as
embedding vector $ \mathbf{W} = (J_1(f_1), \ldots, J_d (f_d)) ,$ we check that
\begin{eqnarray*}
\lefteqn{ E( J_n'(f_n)-J_n(f_n) | \mathbf{W} ) } \notag
\\ &=& - \frac{1}{d} \sum_{i=1}^d \sum_{1\leq i_1<...<i_n\leq d }  \mathbf{1}_{\{ i_1, \ldots, i_n\}}(i)\,\, n! f_n
(i_1 , ..., i_n) E( X_{i_1}\cdot \cdot \cdot X_{i_n} | \mathbf{W}) \notag \\
&=&-\frac{n}d\,J_n(f_n).\label{dollardollar}
\end{eqnarray*}
Thus, with $ \mathbf{W'} = (J_1'(f_1), \ldots, J_d' (f_d)) ,$ the condition \eqref{generalexch} is satisfied, with $\Lambda= (\lambda_{i,j})_{1\leq i,j\leq d}$ being zero
off the diagonal and $\lambda_{n,n} = \frac{n}{d}$ for $n=1, \ldots, d$. Our
 random variable of interest  $F$  is a
linear combination of the elements in $\mathbf{W}$, and we obtain the simple expression
\begin{eqnarray} \label{Lconnection}
E( F'-F |  \mathbf{W}   ) = \frac{1}{d} LF = - \frac{1}{d}  \delta D F.
\end{eqnarray}

This coupling helps to assess the distance to the normal distribution, as follows.

\begin{thm}\label{exchprop}
Let $F\in{\rm dom} D$ be centered and such that $E(F^2) = 1$. Let  $h  \in
{\mathcal C}_b^1$,
 let $Z\sim\mathscr{N}(0,1)$, and let
$(F,F')$ form the exchangeable pair constructed as above, which satisfies \eqref{Lconnection}. Denote by $L'$ the Ornstein-Uhlenbeck operator for the exchanged Rademacher
sequence ${\mathbf{X}}'_d$, and denote by $(L')^{-1}$ its inverse. Then,
\begin{eqnarray*}
|E[h(F)] - E[h(Z)]|&\leq&   4 \| h \|_\infty \sqrt{ {\rm Var} \left[ \frac{d}{2} E\left( (F'-F) \times \big (  (L')^{-1}F'-  L^{-1}F\big) \big | \mathbf{W}\right) \right]
} \nonumber  \\ &&+ \frac{d}{2}   \| h' \|_\infty E\left[ (F' - F)^2\times | (L')^{-1}F'-  L^{-1}F  |\right] .
\end{eqnarray*}
\end{thm}
{\it Proof}. Let $h\in\mathcal{C}^1_b$, and let $g$ denote the solution of the Stein equation \eqref{SteinGaussEq} for $h$. Using antisymmetry, we
have that for all smooth $g$,
\begin{eqnarray}
0 &=&    E \left[ \big(g(F) + g(F')\big)\times\big((L')^{-1}F'-  L^{-1}F\big)\right]  \nonumber  \\ &=&   2E \left[ g(F)\big((L')^{-1}F'-  L^{-1}F\big)\right] +
E\left[\big( g(F') - g(F)\big) \times\big((L')^{-1}F'-  L^{-1}F\big)\right] \nonumber \\ &=& \frac{2}{d} E\big[ F g(F)\big] + E\left[\big( g(F') - g(F) \big)\times
\big((L')^{-1}F'-  L^{-1}F\big)\right];   \label{ivantrick}
\end{eqnarray}
the last equality coming from (\ref{dollardollar}) as well as the definitions of $(L')^{-1}$ and $L^{-1}$. Hence
\begin{eqnarray*}
E\big[g'(F) - F g(F)\big] &=&  E\left[ g'(F) \left(1+  \frac{d}{2}
 E\big[ (F'-F) \times  \big((L')^{-1}F'-  L^{-1}F\big) \big| \mathbf{W}\big]\right) \right] + R,
\end{eqnarray*}
where by Taylor expansion
\begin{eqnarray*}
| R | \leq \frac{d}{4}   \| g''\|_\infty E\big[ (F' - F)^2\times | (L')^{-1}F'-  L^{-1}F |\big] .
\end{eqnarray*}
From \eqref{ivantrick} with $g(x)=x$ we obtain
\begin{eqnarray*}
E \big[ E\big( (F'-F)\times  \big((L')^{-1}F'-  L^{-1}F\big)  \big| \mathbf{W} )\big] &=& E \big[( F' - F )\times \big((L')^{-1}F'-  L^{-1}F\big)\big]  = -\frac{2}{d}
\end{eqnarray*}
by our assumption that $E F^2 = 1$. We now only need to apply the Cauchy-Schwarz inequality and use  the bounds from Lemma \ref{Stein_Lemma_Gauss} to complete  the proof.
\qed

\bigskip
Equation \eqref{generalexch} is key to answering the question on how to construct an exchangeable pair in the case of a random variable
 which admits a finite Rademacher chaos
decomposition. Indeed, one can proceed as follows. Firstly, use $ \mathbf{W} = (J_1(f_1), \ldots, J_d (f_d)) $, where the components come from the chaos decomposition.
Then the results from \cite{ReinertRoellin}  can be applied to assess the distance of $ \mathbf{W} $ to a normal distribution with the same covariance matrix. Note that
due to the orthogonality of the integrals, the covariance matrix will be zero off the diagonal.

\bigskip Finally, note that  formally, using \eqref{Lconnection} and that $L F = \lim_{t \rightarrow 0} \frac{P_t F - F}{t}$ ,
\begin{eqnarray*}
E (F ({\mathbf{X}'}_d) - F ({\mathbf{X}}_d) |  \mathbf{W})
&=&\frac{1}{d}  \lim_{t \rightarrow 0} \frac{P_t F ({\mathbf{X}}_d) - F ({\mathbf{X}}_d)}{t}
= \frac{1}{d}  \lim_{t \rightarrow 0} \frac{1}{t}  [ E (F (\mathbf{X}^t_d)  | {\mathbf{X}}_d ) - F({\mathbf{X}}_d)],
\end{eqnarray*}
where the notation is the same as in Remark \ref{Mehler}. In this sense, our exchangeable pair can be viewed as the limit as $t \rightarrow 0$  of the construction in Remark \ref{Mehler}.

\subsection{From twice differentiable functions to the Wasserstein distance} \label{newsectionbyGesine}

Although many of our results are phrased in terms of smooth test functions, we can also obtain results in Wasserstein distance by mimicking
 the smoothing in \cite{Rinott1996} used for Kolmogorov distance,
 but now apply it to Lipschitz functions. Recall that the Wasserstein distance between the laws of $Y$ and $Z$ is given by
  $$d_W(Y, Z)  = \sup_{h \in {\rm Lip}(1)} \left\vert E[h(Y)]-E[h(Z)] \right\vert, $$
 where ${\rm Lip}(1)$ is the class of all Lipschitz-continuous functions with Lipschitz constant less or equal to 1. Rademacher's Theorem states that a function which is
 Lipschitz continuous
on the real line is Lebesgue-almost everywhere differentiable. For any $h \in {\rm Lip}(1)$ we denote by $h'$ its derivative, which exists almost everywhere.

\begin{cor}
Let $Z\sim\mathscr{N}(0,1)$ and let $F\in{\rm dom} D$ be centered.  Suppose that \eqref{MainBound1} holds for every function $h  \in
{\mathcal C}_b^2$
and that $4( B_1  + B_2  ) \leq 5.$
Then
 \begin{eqnarray*}
 d_W(F, Z) &\leq& \sqrt{2  ( B_1 + B_2 ) (5 + E\vert F \vert )  } .
 \end{eqnarray*}
\end{cor}
{\it Proof.}
 Let $h \in {\rm Lip}(1)$ and, for $t > 0$, define
$$h_t(x) = \int_{-\infty}^{\infty} h( \sqrt{t} y + \sqrt{1-t} x) \phi (y) dy,$$
where $\phi$ denotes the standard normal density. Then we may differentiate and integrate by parts, using that $\phi'(x) = - x \phi(x)$,  to
get
\begin{eqnarray*}
h_t''(x)  &=& \frac{1-t}{\sqrt{t}}
 \int_{-\infty}^{\infty}   y h'( \sqrt{t} y  + \sqrt{1-t} x) \phi \left( y  \right)  dy .
\end{eqnarray*}
Hence for $0 < t < 1$  we may bound
\begin{eqnarray}\label{h2bound}
\| h_t'' \|_\infty &\leq& \frac{1-t}{\sqrt{t}} \| h' \|_\infty \int_{-\infty}^{\infty} \vert y \vert \phi(y) dy \leq \frac{1}{\sqrt{t}}.
\end{eqnarray}
Taylor expansion  gives that for $0 < t
\leq
\frac{1}{2}$ so that $\sqrt{t} \leq \sqrt{1-t}$,
\begin{eqnarray*}
\vert \esp [h(F)] - \esp [h_t(F)]\vert
&
\leq& \left\vert \esp \int  \left\{ h( \sqrt{t} y + \sqrt{1-t} F) - h(\sqrt{1-t} F) \right\} \phi (y) dy \right\vert \\ && +  \esp  \left\vert  h( \sqrt{1-t} F) - h( F)
\right\vert
\\ & \leq & \| h' \|_\infty \sqrt{t} \int | y | \phi(y) dy + \| h' \|_\infty \frac{t}{2 \sqrt{1-t}} \esp |F| \leq \sqrt{t}
\left\{  1 +  \frac{1}{2}  \esp \vert F \vert \right\}.
\end{eqnarray*}
Here we used that $\| h' \|_\infty \leq 1$ and that for $0 < \theta < 1$, we have $(\sqrt{1 - \theta t})^{-1} <  (\sqrt{1 - t})^{-1}$. Similarly,
$
\vert \esp h(Z) - \esp h_t(Z)\vert \leq \frac{3}{2} \sqrt{t}  .
$
Using \eqref{MainBound1} with \eqref{h2bound} together with the triangle inequality
we have
for all   $h \in {\rm Lip}(1)$
\begin{eqnarray*}
 \vert E[h(F)] - E[ h(Z)] \vert
 &\leq &  \frac{1}{\sqrt{t}}
( B_1 (1) + B_2 (1) )
+ \frac{1}{2} \sqrt{t}   \left(  5 + E \vert F \vert \right).
\end{eqnarray*}
Minimising in $t$ gives that the optimal is achieved for
$t = 2(B_1(1) + B_2(1))/(5 + E \vert F \vert)$,
and using this $t$ yields the assertion.
 \qed

\bigskip
To illustrate the result, for Example \ref{ex3.4} we obtain $d_W( F, Z) \leq  \frac{9}{ \sqrt{n}}$
(for $n\geq 6$)
which is of the expected Berry-Ess\'{e}en order.

\section{Normal approximation on a fixed chaos}\label{S : FixedWienerBounds}
\subsection{Explicit upper bounds and CLTs}
We now focus on random variables of the type $F = J_q(f)$, where $J_q(f)$ is the multiple integral defined in Formula (\ref{MWII}), $q\geq 2$ and $f\in\ell^2_0
(\mathbb{N})^{\circ q}$. Due to the definition of the derivative operator $D$, as given in Section \ref{SS : discreteMall}, we know that $F \in {\rm dom}D$. Moreover, by
combining Theorem \ref{firstprop} with formulae (\ref{boundINT1})--(\ref{boundINT2}) one infers that, for every $h\in\mathcal{C}^2_b$  and
for $Z\sim\mathscr{N}(0,1)$,
\begin{eqnarray}
&&\big|E[h(F)] - E[h(Z)]\big| \notag  \\ && \leq \min(4\|h\|_\infty,\|h''\|_{\infty}) \sqrt{E\left[\left(1- \frac1q \| DF\|^2_{\ell^2(\mathbb{N})}\right)^2\right]}
+\frac{20}{3q}\|h''\|_\infty \,E\|DF\|^4_{\ell^4(\N)}.\label{multipleIMPLICIT1}
\end{eqnarray}

The following statement provides an explicit bound of the two expectations appearing in (\ref{multipleIMPLICIT1}).

\begin{thm}\label{MWb}
Fix $q\geq 2$, let $f \in \ell^2_0(\mathbb{N}) ^{\circ q} $, and set $F = J_q (f)$. Then
\begin{eqnarray}
&&E\left\{\left(1-\frac{1}{q} \|DF\|_{\ell^2(\mathbb{N})}^2\right)^2 \right \} = \left |1-q!\|f\|_{\ell^2(\mathbb{N}) ^{\otimes q} }^2\right |^2 \notag \\  &+& q^2
\sum_{p=1}^{q-1}\left \{(p-1)!\binom{q-1}{p-1}^2\right \}^2 (2q-2p)! \,\|\widetilde{f\star_p^p f} \times \mathbf{1}_{\Delta_{2(q-p)}} \|^2_{\ell^2(\mathbb{N}) ^{\otimes
2(q-p)} }  \label{Mww}\\ &\leq& \left |1-q!\|f\|_{\ell^2(\mathbb{N}) ^{\otimes q} }^2\right |^2  + q^2 \sum_{p=1}^{q-1}\left \{(p-1)!\binom{q-1}{p-1}^2\right
\}^2 (2q-2p)! \,\|f\star_p^p f \times \mathbf{1}_{\Delta_{2(q-p)}} \|^2_{\ell^2(\mathbb{N}) ^{\otimes 2(q-p)} }, \notag\\
\label{Mww2}
\end{eqnarray}
and
\begin{eqnarray}
E\|DF\|^4_{\ell^4(\N)} & \leq & q^4 \sum_{p=1}^q \left\{ (p-1)! \binom{q-1}{p-1}^2
 \right\}^2 (2q-2p)!\| f\star_p^{p-1} f \|_{\ell^2(\mathbb{N})^{\otimes (2(q-p)+1)}}
 ^2 \label{Mww4}.
\end{eqnarray}
\end{thm}


\begin{rem}{\rm According to Lemma \ref{L : estimates}, one has that
$$
\| f\star_1^{0} f \|_{\ell^2(\mathbb{N})^{\otimes (2q-1)}}^2 = \| f\star_q^{q-1} f \|_{\ell^2(\mathbb{N})}^2 \leq \| f\|_{\ell^2(\mathbb{N})^{\otimes q}}^2 \times
\max_{j}\sum_{b_1,...,b_{q-1}}f^2(j,b_1,...,b_{q-1}),
$$
and also, by combining this relation with (\ref{MODALX}), for $p=2,...,n-1$,
$$
\| f\star_p^{p-1} f \|_{\ell^2(\mathbb{N})^{\otimes (2(q-p)+1)}}
 ^2 \leq \| f\|_{\ell^2(\mathbb{N})^{\otimes q}}^3 \sqrt{\max_{j}\sum_{b_1,...,b_{q-1}}f^2(j,b_1,...,b_{q-1})}.
$$
These estimates imply that, once $\| f\|_{\ell^2(\mathbb{N})^{\otimes q}}$ is given, the bound on the RHS of (\ref{Mww4}) can be assessed by means of a uniform control on
the ``influence indices'' $\sum_{b_1,...,b_{q-1}}f^2(j,b_1,...,b_{q-1})$, $j\geq 1$. As shown in \cite[Theorem 2.1 and Section 3]{MOO}, when $F=J_q(f)$ depends uniquely
on a finite number of components of $X$, these influence indices can be used to directly measure the distance between the law of $F$ and the law of the random variable,
say $F_G$, obtained by replacing the Rademacher sequence with a i.i.d. Gaussian one. This result roughly implies that the two components of the bound
(\ref{multipleIMPLICIT1}) have a different nature, namely: (\ref{Mww4}) controls the distance between $F$ and its Gaussian-based counterpart $F_G$, whereas
(\ref{Mww})-(\ref{Mww2}) assess the distance between $F_G$ and a standard Gaussian random variable. }
\end{rem}
{\it Proof of Theorem \ref{MWb}}. Observe that, by a classical approximation argument and by virtue of Lemma \ref{L : convLemma}, it is sufficient to consider kernels $f$
with support in a set of the type $\{1,...,N\}^q$, where $N<\infty$ (this is for convenience only, since it allows to freely apply some Fubini arguments). Since $D_jF =
qJ_{q-1}(f(j,\cdot))$, one has that (due to the multiplication formula (\ref{TheProductFormula}))
\begin{eqnarray}
\notag (D_j F) ^2 &=& q^2 \sum_{r=0} ^{q-1} r! \binom{q-1}{r} ^2 J_{2(q-1-r)} \left(\widetilde{f(j,\cdot)\star_r^r f(j,\cdot)}\times {\bf 1}_{\Delta_{2(q-1-r)}}\right)\\
&=&q^2 \sum_{p=1} ^{q} (p-1)! \binom{q-1}{p-1} ^2 J_{2(q-p)} \left(\widetilde{f(j,\cdot)\star_{p-1}^{p-1} f(j,\cdot)}\times {\bf 1}_{\Delta_{2(q-p)}}\right).
\label{Dcomp}
\end{eqnarray}
It follows that, by a Fubini argument,
\begin{eqnarray*}
 \frac1q\|DF\|_{\ell^2(\mathbb{N})}^2 &=&
q \sum_{p=1} ^{q} (p-1)! \binom{q-1}{p-1} ^2 J_{2(q-p)} \left(\widetilde{f\star_p^p f}      \times {\bf 1}_{\Delta_{2(q-p)}}      \right ) \\ &=&
q!\|f\|^2_{\ell^2(\mathbb{N})^{\otimes q}} + q \sum_{p=1} ^{q-1} (p-1)! \binom{q-1}{p-1} ^2 J_{2(q-p)} \left(\widetilde{f\star_p^p f}      \times {\bf
1}_{\Delta_{2(q-p)}}      \right ),
\end{eqnarray*}
and the equality (\ref{Mww}) is obtained by means of orthogonality and isometric properties of multiple integrals. Inequality (\ref{Mww2}) is a consequence of the fact
that, for any kernel $h$, $\|\widetilde{h}\|\leq\|h\|$. To prove (\ref{Mww4}), use again (\ref{Dcomp}) in order to write
\begin{eqnarray*}
E[(D_jF)^4] &=& q^4 \sum_{p=1}^{q} \left\{(p-1)!\binom{q-1}{p-1}^2 \right\}^2 (2q-2p)! \left\|\widetilde{f(j,\cdot)\star_{p-1}^{p-1} f(j,\cdot)}\times {\bf
1}_{\Delta_{2(q-p)}} \right\|_{\ell^2(\mathbb{N})^{\otimes 2(q-p)}} ^2    \\ &\leq& q^4 \sum_{p=1}^{q} \left\{(p-1)!\binom{q-1}{p-1}^2 \right\}^2 (2q-2p)! \left\|
f(j,\cdot)\star_{p-1}^{p-1} f(j,\cdot) \right\|_{\ell^2(\mathbb{N})^{\otimes 2(q-p)}} ^2  .
\end{eqnarray*}
The conclusion is obtained by using the identity
$$
\sum_{j\in\N} \left\| f(j,\cdot)\star_{p-1}^{p-1} f(j,\cdot) \right\|_{\ell^2(\mathbb{N})^{\otimes 2(q-p)}} ^2 = \left\| f\star_{p}^{p-1} f
\right\|_{\ell^2(\mathbb{N})^{\otimes (2(q-p)+1)}} ^2 .
$$
\qed

We deduce the following result, yielding sufficient conditions for CLTs on a fixed chaos.

\begin{prop}\label{P : SuffCLT} Fix $q\geq 2$, and let $F_k=J_q(f_k)$, $k\geq 1$, be a sequence of multiple integrals such that $f_k\in\ell^2_0(\mathbb{N})^{\circ q}$ and
$E(F_k^2) =q!\|f_k\|_{\ell^2(\mathbb{N}) ^{\otimes q} }\rightarrow 1 $. Then, a sufficient condition in order to have that
\begin{equation}\label{MWIcclltt}
F_k \stackrel{\rm Law}{\rightarrow}Z\sim\mathscr{N}(0,1)
\end{equation}
is that
\begin{equation}
\|f_k\star_r^r f_k\|_{\ell^2(\mathbb{N}) ^{\otimes 2(q-r)}} \rightarrow 0, \quad\quad \forall r=1,...,q-1.\label{eqstareq}
\end{equation}
\end{prop}
{\it Proof}. One only needs to combine Theorem \ref{MWb} and Relation (\ref{multipleIMPLICIT1}) with Point 3 of Lemma \ref{L : estimates}. \qed

\begin{rem} {\rm The content of Proposition \ref{P : SuffCLT} echoes the results proved in
\cite{NP}, where it is shown that, on a fixed Gaussian chaos of order $q$, the convergence to zero of the contractions $\|f_k\star_r^r f_k\|$, $r=1,...,q-1$, is necessary
and sufficient in order to have the CLT (\ref{MWIcclltt}). See also \cite{NouPec_ptrf} for bounds, involving norms of contractions, on the normal approximation of the
elements of a Gaussian chaos. See \cite{PSTU} and \cite{PecTaq_Ber} for analogous results concerning the normal approximation of regular functionals of a Poisson measure.
Finally, observe that we do not know at the present time wether the condition (\ref{eqstareq}) is necessary in order to have the CLT (\ref{MWIcclltt}). }
\end{rem}

\subsection{More on finite quadratic forms}\label{SS : quadratic}

When specialized to the case of normalized double integrals, Theorem \ref{firstprop}, Proposition \ref{P : SuffCLT} and Remark \ref{after L : estimates} yield the
following result.

\begin{cor}\label{C : seq2ble}
Let $F_k = J_2(f_k)$, $k\geq 1$, be a sequence of double integrals such that $E(F_k^2)=2\|f_k\|^2_{\ell^2(\mathbb{N})^{\otimes 2}}=1$ (for simplicity). Let $Z\sim
\mathscr{N}(0,1)$. Then, for every $k\geq 1$ and every
$h  \in
{\mathcal C}_b^2$
it holds that
\begin{eqnarray}\label{double bound}
&& \big|E[h(F_k)]-E[h(Z)]\big| \\ && \notag \leq 4 \sqrt{2}\min(4 \|h\|_{\infty},\|h''\|_\infty)\times \|f_k\star_1^1 f_k \times \mathbf{1}_{\Delta_{2}}
\|_{\ell^2(\mathbb{N}) ^{\otimes 2}} + 160 \|h''\|_{\infty}\times\|f_k\star_2^1 f_k \|^2_{\ell^2(\mathbb{N})}\\ && \notag \leq 4 \sqrt{2} \min(4
\|h\|_{\infty},\|h''\|_\infty)\times \|f_k\star_1^1 f_k \times \mathbf{1}_{\Delta_{2}} \|_{\ell^2(\mathbb{N}) ^{\otimes 2}} + 160 \|h''\|_{\infty}\times\|f_k\star_1^1
f_k\times {\bf 1}_{\Delta_2^c} \|^2_{\ell^2(\mathbb{N})^{\otimes 2}},
\end{eqnarray}
and a sufficient condition in order to have the CLT (\ref{MWIcclltt}) is that $\|f_k\star_1^1 f_k \|_{\ell^2(\mathbb{N}) ^{\otimes 2}} \rightarrow 0$.
\end{cor}

Now consider kernels $f_k\in \ell^2_0(\mathbb{N})^{\circ 2}$, $k\geq 1$, such that for every $k$ the following holds: (i) $2\|f_k\|^2_{\ell^2(\mathbb{N})^{\otimes 2}}$ $ =$
$1$, and (ii) the support of $f_k$ is contained in the set $\{1,...,k\}^2$. Then, the random variables
\begin{equation}\label{quu}
F_k =J_2(f_k) = \sum_{1\leq i, j \leq k}f_k(i,j)X_i X_j, \quad k\geq 1,
\end{equation}
are quadratic functionals (with no diagonal terms) of the vectors $(X_1,...,X_k)$. Limit theorems involving sequences such as (\ref{quu}), for general vectors of i.i.d.
random variables, have been intensively studied in the probabilistic literature -- see e.g. \cite{deJongquad}, \cite{GoeTiJTP2002} and the references therein. The
following result, providing a complete characterization of CLTs for sequences such as (\ref{quu}), can be deduced from the main findings contained in \cite{deJongquad}
(but see also \cite{deJongMulti} for extensions to general multilinear functionals).
\begin{prop}\label{P : quadejong}
Let $X$ be the Rademacher sequence considered in this paper, and let $F_k$ be given by (\ref{quu}). For every $k\geq 1$, write $[f_k]$ for the $k\times k$ square
symmetric matrix $\{f_k(i,j) : 1\leq i,j \leq k\}$. Then, the following three conditions are equivalent as $k\rightarrow \infty$:
\begin{itemize}
\item[\rm (a)] $F_k \stackrel{\rm Law}{\rightarrow} Z\sim \mathscr{N}(0,1)$; \item[\rm (b)] ${\bf Trace}\{[f_k]^4\}= \|f_k\star_1^1 f_k
    \|^2_{\ell^2(\mathbb{N})^{\otimes 2}}\rightarrow 0$; \item[\rm (c)] $E(F_k^4) \rightarrow E(Z^4) = 3$.
\end{itemize}
\end{prop}

\begin{rem}{\rm Condition (b) in the previous statement is often replaced by the following:
$$ {\rm (b')} \,\,\, {\bf Trace}\{[f_k]^4\}\rightarrow 0 \quad \text{and} \quad \max_{j\leq k}\sum_{i=1}^k f_k(i,j)^2\rightarrow 0.$$ However,
Relation (\ref{traceequality}) shows that the second requirement in (${\rm b}'$) is indeed redundant. }
\end{rem}

In \cite{Chatterjee_AOP}, Chatterjee proves that, in Wasserstein distance,
\begin{equation}\label{chatterbound}
d_W(F_k,Z)\leq \sqrt{\frac12 {\bf Trace}\{[f_k]^4\} }+ \frac52\sum_{j=1}^k \big[\sum_{i=1}^k f_k(i,j)\big]^{\frac32}.
\end{equation}
Note that
\begin{equation}\label{chatterbound2}
\max_{j\leq k} \big[\sum_{i=1}^k f_k(i,j)^2\big]^{\frac32} \leq \sum_{j=1}^k \big[\sum_{i=1}^k f_k(i,j)^2\big]^{\frac32}\leq \frac12 \max_{j\leq k}\big[\sum_{i=1}^k
f_k(i,j)^2\big]^{\frac12},
\end{equation}
and due e.g. to (\ref{traceequality}), Relation (\ref{chatterbound}) gives an alternate proof of the implication (b) $\rightarrow$ (a) in Proposition \ref{P : quadejong}.
Another proof of the same implication can be deduced from the following result, which is a direct consequence of Corollary \ref{C : seq2ble} and Lemma \ref{L : estimates}
(Point 3).

\begin{cor} Let the notation and assumptions of this section prevail. Then, for every $k$ and for every
 $h  \in
{\mathcal C}_b^2$,
one has
that
\begin{eqnarray*}
&& |E[h(F_k)]-E[h(Z)]| \\ && \notag \leq 4 \sqrt{2}\min(4 \|h\|_{\infty},\|h''\|_\infty) \|f_k\star_1^1 f_k \times \mathbf{1}_{\Delta_{2}} \|_{\ell^2(\mathbb{N}) ^{\otimes
2}} + 160 \|h''\|_{\infty}\|f_k\star_1^1 f_k \times \mathbf{1}_{\Delta_{2}^c}\|^2_{\ell^2(\mathbb{N})^{\otimes 2}}\\ &&
\leq 4\sqrt{2}\min(4 \|h\|_{\infty},\|h''\|_\infty)
\sqrt{{\bf Trace}\{[f_k]^4\}}  + 160 \|h''\|_{\infty} {\bf Trace}\{[f_k]^4\}.
\end{eqnarray*}
\end{cor}

\section{Sums of single and double integrals and applications to weighted runs}\label{S : 2runs}
\subsection{Motivation: infinite weighted runs}
Due to their simple dependence structure, runs lend themselves as good test examples for normal approximations; for applications  see for example \cite{Balakrishnan2002}.
The first Berry-Ess\'{e}en bound for overlapping success runs was derived in \cite{Godbole}.
 A multivariate normal approximation for the count runs of finite length  in the context of Stein's method can be found for example
in \cite{ReinertRoellin}. Weighted runs can be viewed as a special case of weighted U-statistics; a normal approximation using Stein's method is available for example in
\cite{Rinott1997}. Typically these studies focus on counting runs in a finite window of a possibly infinite sequence, see for example \cite{Balakrishnan2002} for details.
In contrast, our method allows to consider infinite weighted sums of runs.

\smallskip

Let $\xi=\{\xi_n:\,n\in \Z\}$ be a standard {\sl Bernoulli sequence}. This means that $\xi$ is a collection of i.i.d. random
variables such that $P[\xi_1=1]=P[\xi_1=0]=1/2$. In the sequel, we will suppose without loss of generality that $\xi=\frac12(1+X)$
where $X=\{X_n:\,n\in\Z\}$ is a two-sided {\sl Rademacher} sequence, obtained e.g. by juxtaposing two independent
Rademacher sequences of the kind introduced in Section \ref{setup}. Note that we chose $\mathbb{Z}$ as a parameter set,
in order to simplify the forthcoming discussion. Also, in what follows we will implicitly use the fact that all the
results and bounds established in the previous sections extend to this setting, by simply replacing sums over $\mathbb{N}$ with
 sums over $\mathbb{Z}$ (the proof of this elementary fact is left to the reader).
Now fix an integer $d\geq 1$ as well as a sequence $\{\alpha^{(n)}:\,n\geq 1\}$ of elements of $\ell^2(\Z)$. In this section,
we  study the normal approximation of the sequence $\{G_n:\,n\geq 1\}$ defined by
\begin{equation}\label{Gn}
G_n=\sum_{i\in\Z} \alpha^{(n)}_i \xi_i \ldots \xi_{i+d}.
\end{equation}
We call $G_n$ an {\sl infinite weighted $d$-run}. Actually, we will rather focus on the normalized version of $G_n$, namely
\begin{equation}\label{Fn}
F_n=\frac{G_n-E(G_n)}{\sqrt{{\rm Var}G_n}}.
\end{equation}
Observe that, in order to apply Theorem \ref{firstprop}, one needs to find the chaotic expansion of $F_n$ or, equivalently, that of $G_n$. Using the identity
$\xi=\frac12(1+X)$, it is immediately deduced, see also Section \ref{SS : how?}-(ii), that
\begin{eqnarray}
G_n&=&2^{-(d+1)}\sum_{i\in\Z} \alpha_i^{(n)} (1+X_i)\ldots(1+X_{i+d})\notag\\ &=&2^{-(d+1)}\sum_{i\in\Z}\alpha_i^{(n)} \sum_{I\subset\{i,\ldots,i+d\}} X_{i_1}\ldots
X_{i_{|I|}}\quad\mbox{(with $I=\{i_1,\ldots,i_{|I|}\}$)}\notag\\ &=&2^{-(d+1)} \sum_{r=0}^{d+1} \sum_{i\in\Z}\alpha_i^{(n)}
\sum_{\substack{I\subset\{i,\ldots,i+d\}\\|I|=r}} X_{i_1}\ldots X_{i_r}\notag\\ &=&E(G_n)+ \sum_{r=1}^{d+1} J_r\left(2^{-(d+1)}\sum_{i\in\Z}\alpha_i^{(n)}
\sum_{\substack{I\subset\{i,\ldots,i+d\}\\|I|=r}} \widetilde{{\bf 1}_{\{i_1\}}\otimes \ldots \otimes{\bf 1}_{\{i_r\}}}\right),\label{deco}
\end{eqnarray}
where the tilde indicates a symmetrization. In particular, it is now immediate to compute $DG_n$ and $DL^{-1}G_n$, by using the definitions given in Section \ref{SS :
discreteMall}. However, since the analysis of (\ref{Gn}) is meant to be only an illustration, from now on we will focus on the case where $d=1$ ($2$-runs). The general
case could be handled in a similar way (at the cost of a quite cumbersome notation).

\subsection{Normal approximation of sums of single and double integrals}
We will deduce a bound for the quantity $\big| E[h(F_n)]-E[h(Z)]\big|$ (where $Z\sim\mathscr{N}(0,1)$ and $F_n$ is given by (\ref{Fn})) from the following result, which
can be seen as a particular case of Theorem \ref{firstprop}.
\begin{prop}\label{sum1+2}
{\rm (Sum of a single and a double integral)} Let $F=J_1(f)+J_2(g)$ with $f\in\ell^2(\Z)$ and $g\in\ell_0^2(\Z)^{\circ 2}$.
Assume that $\sum_{i\in\Z}
\big|g(i,k)\big|<\infty$ for all $k\in\Z$. Also, suppose (for simplicity) that ${\rm Var} F=1$,
and let  $h  \in
{\mathcal C}_b^2$.
Then
\begin{eqnarray}
\big|E[h(F)]-E[h(Z)]\big|&\leq& \min(4\|h\|_\infty,\|h''\|_\infty) \big( 2\sqrt{2}\|g\star_1^1g{\bf 1}_{\Delta_2}\|_{\ell^2(\Z)^{\otimes 2}}
+3\|f\star_1^1g\|_{\ell^2(\Z)} \big)\notag\\ &&+\frac{160}{3}\|h''\|_\infty\sum_{k\in\Z} \left[f^4(k) + 16\left(\sum_{i\in\Z} |g(i,k)|\right)^4\right].\label{rhsjs}
\end{eqnarray}
\end{prop}
\begin{rem}{\rm By applying successively Fubini and Cauchy-Schwarz theorems, one has that
\begin{equation}\label{join}
\|f\star_1^1g\|^2_{\ell^2(\Z)}
=\sum_{i,j\in\Z}f(i)f(j)\,\,g\star_1^1g(i,j)
\leq \|f\|^2_{\ell^2(\mathbb{Z})}\|g\star_1^1 g\|_{\ell^2(\mathbb{Z})^{\otimes 2}}.
\end{equation}
This inequality has an interesting consequence. Suppose indeed that the sequences $J_1(f_n)$ and $J_2(g_n)$, $n\geq 1$, are such that, as $n\rightarrow\infty$: (a)
$\|f_n\|_{\ell^2(\mathbb{Z})}\rightarrow 1$, (b) $2\|g_n\|^2_{\ell^2(\mathbb{Z})^{\otimes 2}}\rightarrow 1$, (c) $\sum_{k\in\Z} f_n^4(k)\rightarrow 0$, (d) $\sum_{k\in\Z}
\left(\sum_{i\in\Z} |g_n(i,k)|\right)^4 \rightarrow 0$, and
 (e) $\|g_n\star_1^1 g_n\|_{\ell^2(\mathbb{Z})^{\otimes 2}}\rightarrow 0$. Then the estimate (\ref{join}) and Proposition \ref{sum1+2} imply that, for every
 $(\alpha,\beta)\neq(0,0)$, the sequence \begin{equation*}
\frac{\alpha J_1(f_n) +\beta J_2(g_n)}{\sqrt{\alpha^2 +\beta^2}}, \quad n\geq 1,
\end{equation*}
converges in law to $Z\sim\mathscr{N}(0,1)$, and therefore that the vectors $(J_1(f_n), J_2(g_n))$, $n\geq 1$, jointly converge in law towards a two-dimensional centered
i.i.d. Gaussian vector with
unit variances. Note that each one of conditions (a)--(e) involves {\sl separately} one of the two kernels $f_n$ and $g_n$. See
\cite{NouPeRev} and \cite[Section 6]{PSTU}, respectively, for several related results in a Gaussian and in a Poisson framework. }
\end{rem}

\noindent{\it Proof of Proposition \ref{sum1+2}}. Firstly, observe that $D_kF= 2J_1\big(g(\cdot,k)\big)+f(k)=2\sum_{i\in\Z} g(i,k)X_i+f(k)$ so that, using $(a+b)^4\leq
8(a^4+b^4)$,
$$
E\big|D_k F|^4\leq 128 E \left(\sum_{i\in\Z} g(i,k)X_i\right)^4 + 8f^4(k)\leq 128 \left(\sum_{i\in\Z} \big|g(i,k)\big|\right)^4 + 8f^4(k)<\infty.
$$
We are thus left to bound
$B_1$ and $B_2$
in Theorem \ref{firstprop}, taking into account the particular form of $F$. We have $-L^{-1}F=\frac12J_2(g)+J_1(f)$ so
that $-D_kL^{-1}F= J_1\big(g(\cdot,k)\big)+f(k)$. Consequently,  with  the multiplication formula (\ref{TheProductFormula}), we get
\begin{eqnarray*}
&&\langle DF,-DL^{-1}F\rangle_{\ell^2(\Z)}\\ &=&2\sum_{k\in\Z} J_1\big(g(\cdot,k)\big)^2+3\sum_{k\in\Z}f(k)J_1\big(g(\cdot,k)\big) + \|f\|^2_{\ell^2(\Z)}\\
&=&2\sum_{k\in\Z} J_2\big(g(\cdot,k)\otimes g(\cdot,k){\bf 1}_{\Delta_2}\big) +3\sum_{k\in\Z}f(k)J_1\big(g(\cdot,k)\big) + \|f\|^2_{\ell^2(\Z)} +
2\|g\|^2_{\ell^2(\Z)^{\otimes 2}}\\ &=&2 J_2(g\star_1^1g{\bf 1}_{\Delta_2}) +3J_1(f\star_1^1g) + \|f\|^2_{\ell^2(\Z)} + 2\|g\|^2_{\ell^2(\Z)^{\otimes 2}}
\end{eqnarray*}
so that
\begin{eqnarray*}
E\big|\langle DF,-DL^{-1}F\rangle_{\ell^2(\Z)}-{\rm Var}F\big|^2 &=&E\left|2 J_2(g\star_1^1g{\bf 1}_{\Delta_2}) +3J_1(f\star_1^1g) \right|^2 \\ &=&8\|g\star_1^1g{\bf
1}_{\Delta_2}\|^2_{\ell^2(\Z)^{\otimes 2}} +9\|f\star_1^1g\|_{\ell^2(\Z)}^2.
\end{eqnarray*}
Hence,
\begin{eqnarray*}
B_1&\leq&
\sqrt{8\|g\star_1^1g{\bf 1}_{\Delta_2}\|^2_{\ell^2(\Z)^{\otimes 2}} +9\|f\star_1^1g\|_{\ell^2(\Z)}^2}\\ &\leq&
 \big( 2\sqrt{2}\|g\star_1^1g{\bf 1}_{\Delta_2}\|_{\ell^2(\Z)^{\otimes 2}} +3\|f\star_1^1g\|_{\ell^2(\Z)} \big).
\end{eqnarray*}
Now, let us consider
$B_2$. We have
$$
\big|D_kF\big| \leq |f(k)| + 2\sum_{i\in\Z} |g(i,k)|.
$$
Similarly,
\begin{eqnarray*}
\big|D_kL^{-1}F\big| &=& \big|f(k)+J_1\big(g(\cdot,k)\big)\big| =\left|f(k)+\sum_{i\in\Z} g(i,k)X_i\right|\\ &\leq& |f(k)| + \sum_{i\in\Z} |g(i,k)|\leq |f(k)| +
2\sum_{i\in\Z} |g(i,k)|.
\end{eqnarray*}
Still using $(a+b)^4\leq 8(a^4+b^4)$, we deduce
\begin{eqnarray*}
\sum_{k\in\Z} \big|D_kL^{-1}F\big|\times \big|D_kF\big|^3 &\leq& \sum_{k\in\Z} \left( |f(k)| + 2\sum_{i\in\Z} |g(i,k)|\right)^4\\ &\leq& 8 \sum_{k\in\Z} \left[f^4(k) +
16\left(\sum_{i\in\Z} |g(i,k)|\right)^4\right].
\end{eqnarray*}
Hence
$$
B_2 \leq \frac{160}{3}
\sum_{k\in\Z} \left[f^4(k) + 16\left(\sum_{i\in\Z} |g(i,k)|\right)^4\right]
$$
and the desired conclusion follows by applying Theorem \ref{firstprop}.\qed
\subsection{Bounds for infinite $2$-runs}
When $d=1$, Proposition \ref{sum1+2} allows to deduce the following bound for the normal approximation of the random variable $F_n$ defined in (\ref{Fn}).
\begin{prop}
Let $\{F_n:\,n\geq 1\}$ be the sequence defined by $F_n=\frac{G_n-E(G_n)}{\sqrt{{\rm Var}G_n}}$ with
$$
G_n=\sum_{i\in\Z}\alpha_i^{(n)}\xi_i\xi_{i+1}.
$$
Here, $\xi=\{\xi_n:\,n\in \Z\}$ stands for the standard {\it Bernoulli} sequence and $\{\alpha^{(n)}:\,n\geq 1\}$
 is a given sequence of elements of $\ell^2(\Z)$.
Consider a function $h  \in
{\mathcal C}_b^2$.
Then, for $Z\sim\mathscr{N}(0,1)$,
\begin{eqnarray}
\big|E[h(F)]-E[h(Z)]\big|&\leq& \frac{7}{16}\times \frac{\min(4\|h\|_\infty,\|h''\|_\infty) }{{\rm Var}G_n}\times \sqrt{\sum_{i\in\Z}(\alpha_i^{(n)})^4}\label{wysiwyg}
\\
&& +\frac{35}{24}\times \frac{\|h''\|_\infty}{({\rm Var}G_n)^2}\times \sum_{i\in\Z}(\alpha_i^{(n)})^4\notag
\end{eqnarray}
with
\begin{equation}\label{vargn}
{\rm Var}G_n=\frac3{16}\sum_{i\in\Z}(\alpha_i^{(n)})^2 + \frac18\sum_{i\in\Z}\alpha_i^{(n)}\alpha_{i+1}^{(n)}.
\end{equation}
It follows that a sufficient condition to have $F_n \stackrel{\rm Law}{\rightarrow} Z$ is that
$$
\sum_{i\in\Z}(\alpha_i^{(n)})^4=o\left(({\rm Var}G_n)^2\right)\quad \mbox{as } n\to\infty .
$$
\end{prop}
{\it Proof}. Identity (\ref{vargn}) is easily verified. On the other hand, by (\ref{deco}), we have
$$
F_n=\frac{G_n - E(G_n)}{ \sqrt{{\rm Var}G_n} } = J_1(f) + J_2(g),
$$
with
\begin{eqnarray*}
f&=&\frac1{4\sqrt{{\rm Var}G_n}}\sum_{a\in\Z} \alpha_a^{(n)} \big( {\bf 1}_{\{a\}} + {\bf 1}_{\{a+1\}}\big)\\ g&=&\frac1{8\sqrt{{\rm Var}G_n}} \sum_{a\in\Z}
\alpha_a^{(n)} \big( {\bf 1}_{\{a\}}\otimes {\bf 1}_{\{a+1\}} + {\bf 1}_{\{a+1\}}\otimes {\bf 1}_{\{a\}}\big).
\end{eqnarray*}
Now, let us compute each quantity appearing in the RHS of (\ref{rhsjs}). If $i\neq j$ then
\begin{eqnarray*}
(g\star_1^1 g)(i,j)&=& \frac{1}{64{\rm Var}G_n} \sum_{a,b,k\in\Z} \alpha_a^{(n)}\alpha_b^{(n)} \bigg( {\bf 1}_{\{a\}}(i) {\bf 1}_{\{a+1\}}(k) + {\bf 1}_{\{a+1\}}(i) {\bf
1}_{\{a\}}(k)\bigg)\\ &&\hskip3.5cm\times \bigg( {\bf 1}_{\{b\}}(j) {\bf 1}_{\{b+1\}}(k) + {\bf 1}_{\{b+1\}}(j) {\bf 1}_{\{b\}}(k)\bigg)\\ &=&\frac{1}{64{\rm Var}G_n}
\left( \alpha_i^{(n)} \alpha_{i+1}^{(n)} {\bf 1}_{\{ j=i+2\}} + \alpha_j^{(n)} \alpha_{j+1}^{(n)} {\bf 1}_{\{ j=i-2\}} \right).
\end{eqnarray*}
Hence
\begin{eqnarray*}
&&\|g\star_1^1 g{\bf 1}_{\Delta_2}\|_{\ell^2(\Z)^{\otimes 2}}\\ &=&\frac{1}{64{\rm Var}G_n} \sqrt{\sum_{i,j\in\Z}\left[ (\alpha_i^{(n)})^2 (\alpha_{i+1}^{(n)})^2 {\bf
1}_{\{ j=i+2\}} + (\alpha_j^{(n)})^2 (\alpha_{j+1}^{(n)})^2 {\bf 1}_{\{ j=i-2\}}\right]}\\ &=&\frac{\sqrt{2}}{64{\rm Var}G_n} \sqrt{\sum_{i\in\Z} (\alpha_i^{(n)})^2
(\alpha_{i+1}^{(n)})^2} .
\end{eqnarray*}
We have
\begin{eqnarray*}
(f\star_1^1 g)(i)&=&\frac{1}{32{\rm Var}G_n} \sum_{a,b,k\in\Z} \alpha_a^{(n)}\alpha_b^{(n)} \big( {\bf 1}_{\{a\}}(k) + {\bf 1}_{\{a+1\}}(k) \big)\\ &&\hskip3.5cm\times
\bigg( {\bf 1}_{\{b\}}(i) {\bf 1}_{\{b+1\}}(k) + {\bf 1}_{\{b+1\}}(i) {\bf 1}_{\{b\}}(k)\bigg)\\ &=&\frac{1}{32{\rm Var}G_n} \left( \alpha_i^{(n)} \alpha_{i+1}^{(n)} +
(\alpha_{i-1}^{(n)})^2 + (\alpha_{i}^{(n)})^2 + \alpha_{i-1}^{(n)}\alpha_{i-2}^{(n)} \right).
\end{eqnarray*}
Hence, using $(a+b+c+d)^2\leq 4(a^2+b^2+c^2+d^2)$,
\begin{eqnarray*}
&&\|f\star_1^1 g\|_{\ell^2(\Z)}\\ &=&\frac{1}{32{\rm Var}G_n}\sqrt{ \sum_{i\in\Z}\left[ \alpha_i^{(n)} \alpha_{i+1}^{(n)} +  (\alpha_{i-1}^{(n)})^2 + (\alpha_{i}^{(n)})^2
+ \alpha_{i-1}^{(n)}\alpha_{i-2}^{(n)} \right]^2 }\\ &\leq&\frac{1}{16{\rm Var}G_n}\sqrt{ \sum_{i\in\Z}(\alpha_i^{(n)})^2 (\alpha_{i+1}^{(n)})^2 +
\sum_{i\in\Z}(\alpha_{i-1}^{(n)})^4 + \sum_{i\in\Z}(\alpha_{i}^{(n)})^4 + \sum_{i\in\Z}(\alpha_{i-1}^{(n)})^2(\alpha_{i-2}^{(n)})^2 } \\ &=&\frac{\sqrt{2}}{16{\rm
Var}G_n}\sqrt{ \sum_{i\in\Z}(\alpha_i^{(n)})^2 (\alpha_{i+1}^{(n)})^2 + \sum_{i\in\Z}(\alpha_{i}^{(n)})^4 } .
\end{eqnarray*}
We have, using $(a+b)^4\leq 8(a^4+b^4)$,
\begin{eqnarray*}
\sum_{k\in\Z}f^4(k)&=&\frac{1}{256({\rm Var}G_n)^2} \sum_{k\in\Z} \left[ \sum_{a\in\Z}\alpha_a^{(n)} \big( {\bf 1}_{\{a\}}(k) + {\bf 1}_{\{a+1\}}(k)\big)\right]^4\\
&=&\frac{1}{256({\rm Var}G_n)^2} \sum_{k\in\Z} \big( \alpha_k^{(n)} + \alpha_{k-1}^{(n)}\big)^4\\ &\leq&\frac{1}{16({\rm Var}G_n)^2} \sum_{k\in\Z} \big( \alpha_k^{(n)}
\big)^4.
\end{eqnarray*}
Finally, still using $(a+b)^4\leq 8(a^4+b^4)$,
\begin{eqnarray*}
&&\sum_{k\in\Z}\left[\sum_{i\in\Z}|g(i,k)|\right]^4
\\ &\leq&\frac{1}{4096({\rm Var}G_n)^2}\sum_{k\in\Z}
\left[\sum_{i,a\in\Z} \big|\alpha_a^{(n)}\big| \big| {\bf 1}_{\{a\}}(i) {\bf 1}_{\{a+1\}}(k) + {\bf 1}_{\{a+1\}}(i) {\bf 1}_{\{a\}}(k)\big|\right]^4\\
&=&\frac{1}{4096({\rm Var}G_n)^2}\sum_{k\in\Z} \left[\big|\alpha_{k-1}^{(n)}\big|+\big|\alpha_{k}^{(n)}\big|\right]^4\\ &\leq&\frac{1}{256({\rm
Var}G_n)^2}\sum_{k\in\Z}(\alpha_k^{(n)})^4.
\end{eqnarray*}
Now, the desired conclusion follows by plugging all these estimates in (\ref{rhsjs}), after observing that $\sum_{i\in\Z} (\alpha_i^{(n)})^2 (\alpha_{i+1}^{(n)})^2 \leq
\sum_{i\in\Z} (\alpha_i^{(n)})^4$, by the Cauchy-Schwarz inequality. \qed

\begin{example}
{\rm
\begin{enumerate}
\item Choose $\alpha_i^{(n)}={\bf 1}_{\{1,\ldots,n\}}(i)$. Then $\sum_{i\in\Z}(\alpha_i^{(n)})^4 = n$ and $${\rm Var}G_n\geq \frac3{16}\sum_{i\in\Z}(\alpha_i^{(n)})^2
    = \frac{3n}{16},$$ so that (\ref{wysiwyg}) gives a bound of order $n^{-1/2}$ overall. This is the same order as obtained in \cite{ReinertRoellin}, also for smooth
    test functions, while \cite{Rinott1997} obtain a bound of this order in Kolmogorov distance under suitable conditions on the weights. \item Choose
    $\alpha_i^{(n)}=i^{-1}\,{\bf 1}_{\{n,n+1,\ldots\}}(i)$. Then $\sum_{i\in\Z}(\alpha_i^{(n)})^4 = O(n^{-3})$ and $${\rm Var}G_n\geq
    \frac3{16}\sum_{i\in\Z}(\alpha_i^{(n)})^2 \sim_{n\to\infty}\frac3{16n}$$ so that (\ref{wysiwyg})  also gives a bound of order $n^{-1/2}$ overall.
\end{enumerate}
}
\end{example}

\section{Multiple integrals over sparse sets}\label{S : App}
\subsection{General results}\label{SSGenBJ}
Fix $d\geq 2$. Let $F_N$, $N\geq 1$, be a sequence of subsets of $\mathbb{N}^d$ such that the following three properties are satisfied for every $N\geq 1$: (i) $F_N \neq
\emptyset$, (ii) $F_N \subset \Delta_d^N$ (as defined in (\ref{simplex2})), that is, $F_N$ is contained in $\{1,\ldots,N\}^d$ and has no diagonal components, and (iii)
$F_N$ is a symmetric set, in the sense that every $(i_1,...,i_d)\in F_N$ is such that $(i_{\sigma(1)},...,i_{\sigma(d)})\in F_N$ for every permutation $\sigma$ of the set
$\{1,...,d\}$. Let $X$ be the infinite Rademacher sequence considered in this paper. Given sets $F_N$ as at points (i)--(iii), we shall consider the sequence of
multilinear forms
\begin{equation}\label{MF}
\widetilde{S}_N = [d!\times|F_N|]^{-\frac12}\sum_{(i_1,...,i_d) \in F_N} X_{i_1}\cdot\cdot\cdot X_{i_d} = J_d (f_N), \quad N\geq 1,
\end{equation}
where $|F_N|$ stands for the cardinality of $F_N$,
and
$$
f_N(i_1,...,i_d) := [d!\times |F_N|]^{-\frac12}\times {\bf 1}_{F_N} (i_1,...,i_d).
$$
Note that $E(\widetilde{S}_N) =0$ and $E(\widetilde{S}_N^2) =1$ for every $N$. In the paper \cite{BleiJanson}, Blei and Janson studied the problem of finding conditions on the set $F_N$, in
order to have that the CLT
\begin{equation}\label{BJclt}
\widetilde{S}_N \stackrel{\rm Law}{\rightarrow} Z\sim \mathscr{N}(0,1), \quad N\rightarrow \infty,
\end{equation}
holds.

\begin{rem}{\rm Strictly speaking, Blei and Janson use the notation $F_N$ in order to indicate the \textsl{restriction to the simplex}
$\{(i_1,...,i_d):i_1<i_2<...<i_d\}$ of a set verifying Properties (i)--(iii) above. }
\end{rem}

In order to state Blei and Janson's result, we need to introduce some more notation.

\medskip

\noindent{\bf Remark on notation.} In what follows, we will write ${\bf a}_k$
to indicate vectors ${\bf a}_k = (a_1,...,a_k)$ belonging to a set of the type
$\{1,...,N\}^k=:[N]^k$, for some $k,N\geq 1$. We will regard these objects both as vectors and sets, for instance: an expression of the type ${\bf a}_k \cap{\bf i}_l =
\emptyset$, means that the two sets $\{a_1,...,a_k\}$ and $\{i_1,...,i_l\}$ have no elements in common; when writing $j\in {\bf a}_k$, we mean that $j=a_r$ for some
$r=1,...,k$; when writing ${\bf a}_k \subset{\bf i}_l$ ($k\leq l$), we indicate that, for every $r=1,...,k$, one has $a_r =i_s$ for some $s=1,...,l$. When a vector ${\bf
a}_k$ enters in a sum, we will avoid to specify ${\bf a}_k\in [N]^k$, whenever the domain of summation $[N]^k$ is clear from the context.

\medskip

Given $N\geq 1$ and an index $j\in [N]$, we set
$$
F^*_{N,j} = \{{\bf i}_d \in F_N : j\in {\bf i}_d\}.
$$
For every $N$, the set $F^{\#}_N\subset F_N \times F_N$ is defined as the collection of all pairs $({\bf i}_d, {\bf k}_d) \in F_N \times F_N$ such that: (a) ${\bf
i}_d\cap {\bf k}_d = \emptyset$, and (b) there exists $p=1,...,d-1$, as well as ${\bf i}'_p\subset {\bf i}_d $ and ${\bf k}'_p\subset {\bf k}_d $ such that
$$
(({\bf k}'_p,\,{\bf i}_d \setminus {\bf i}'_p),\,({\bf i}'_p,\,{\bf k}_d \setminus {\bf k}'_p)) \in F_N\times F_N,
$$
where ${\bf i}_d \setminus {\bf i}'_p$ represents the element of $[N]^{d-p}$ obtained by eliminating from ${\bf i}_d$ the coordinates belonging to ${\bf i}'_p$, and
$({\bf k}'_p,\,{\bf i}_d \setminus {\bf i}'_p)$ is the element of $[N]^d$ obtained by replacing ${\bf i}'_p$ with ${\bf k}'_p$ in ${\bf i}_d$ (an analogous description
holds for $({\bf i}'_p,\,{\bf k}_d \setminus {\bf k}'_p))$. In other words, the $2d$ indices $i_1,\ldots,i_d,j_1,\ldots,j_d$ can be partitioned in at least two ways into
elements of $F_N$.

\begin{thm}\label{T : Blei and Janson}{\rm (\cite[Th. 1.7]{BleiJanson})} Let the above
notation and assumptions prevail, and suppose that
\begin{eqnarray}
&& \lim_{N\rightarrow\infty} \max_{j\leq N} \frac{|F^*_{N,j}|}{|F_{N}|} = 0, \quad \text{and} \label{BJneg1}\\ && \lim_{N\rightarrow\infty} \frac{|F^\#_{N}|}{|F_{N}|^2} =
0 \label{BJneg2}.
\end{eqnarray}
Then Relation (\ref{BJclt}) holds, with convergence of all moments.
\end{thm}

\begin{rem}{\rm As pointed out in \cite{BleiJanson}, Condition (\ref{BJneg2}) can be described
as a weak ``sparseness condition'' (see e.g. \cite{BleiBook}). See also \cite[Th. 1.7]{BleiJanson} for a converse statement. }
\end{rem}

The principal achievement of this section is the following refinement of Theorem \ref{T : Blei and Janson}.

\begin{thm}\label{T : BerryEsseenBJ} Under the above notation and assumptions, consider a
function $h  \in
{\mathcal C}_b^2$.
Then there exist two universal constants $C_1$ and $C_2$, depending only on $d$, $\|h\|_\infty$ and
$\|h''\|_\infty$, such that, for every $N\geq 1$ and for $Z\sim\mathscr{N}(0,1)$,
\begin{equation}\label{BJbound}
 |E[h(Z)] - E[h(\widetilde{S}_N)]|
 \leq C_1 \frac{|F^\#_{N}|^{\frac12}}{|F_{N}|}
+C_2 \left[\max_{j\leq N}\frac{|F^*_{N,j}|}{|F_{N}|}\right]^{\frac14}.
\end{equation}
\end{thm}
{\it Proof.} By combining (\ref{multipleIMPLICIT1}) with Theorem \ref{MWb}, we know that there exist universal combinatorial constants $\alpha_p^d>0$, $p=1,...,d-1$ and
$\beta_l^d>0$, $l=1,...,d$, such that
\begin{eqnarray}\label{G1}
 \big|E[h(\widetilde{S}_N)]  - E[h(Z)]\big| &\leq &
\|h\|_\infty \sum_{p=1}^{d-1} \alpha_p^d\, \|f_N\star^p_p f_N \, {\bf 1}_{\Delta_{2(d-p)}^N}\|_{\ell^2(\mathbb{N})^{\otimes 2(d-p)}} \\&& + \|h''\|_\infty\sum_{l=1}^{d-1}
\beta_l^d \, \|f_N\star^{l-1}_{l} f_N \|^2_{\ell^2(\mathbb{N})^{\otimes 2(d-l)+1}}. \label{G2}
\end{eqnarray}
We now evaluate each norm appearing in (\ref{G1})--(\ref{G2}). According to the second and the third points of Lemma \ref{L : estimates} and using the fact that
$\|f_N\|^2_{\ell^2(\mathbb{N})^{\otimes d}} = (d!)^{-1}$, one has that
\begin{eqnarray} \notag
\| f_N\star_1^0 f_N \|^2_{\ell^2(\mathbb{N})^{\otimes (2d-1)}}& = &\| f_N\star_d^{d-1} f_N \|^2_{\ell^2(\mathbb{N})} \\ & \leq & \frac{1}{d!}\max_{j\leq N} \sum_{{\bf
b}_{d-1}} f_N(j,{\bf b}_{d-1})^2 =\frac1{d\times d!}\max_{j\leq N} \frac{|F^*_{N,j}|}{|F_{N}|}.\label{FRLA1}
\end{eqnarray}
Also, by combining the two inequalities in the second point of Lemma \ref{L : estimates}, we deduce that, for every $l=2,...,d-1$,
\begin{eqnarray}\label{FRLA2}
\| f_N\star_l^{l-1} f_N \|^2_{\ell^2(\mathbb{N})^{\otimes 2(d-l)+1}}&\leq & \frac{1}{d!}\| f_N\star_d^{d-1} f_N \|_{\ell^2(\mathbb{N})} \\ & \leq &
\frac{1}{(d!)^{3/2}\sqrt{d}} \left[\max_{j\leq N}\frac{|F^*_{N,j}|}{|F_{N}|}\right]^{\frac12} \leq \frac{1}{(d!)^{3/2}\sqrt{d}} \left[\max_{j\leq
N}\frac{|F^*_{N,j}|}{|F_{N}|}\right]^{\frac14}.\notag
\end{eqnarray}
Now fix $p=1,...,d-1$. One has that
\begin{eqnarray}
\label{cool}&& \|f_N\star^p_p f_N \,{\bf 1}_{\Delta_{2(d-p)}^N}\|^2_{\ell^2(\mathbb{N})^{\otimes 2(d-p)}}\\ && = \frac1{(d!|F_N|)^2} {\sum_{\substack{{\bf a}_{d-p},{\bf
b}_{p} \\ {\bf x}_{d-p},{\bf y}_{p}}}}\!\!{\bf 1}_{F_N}({\bf a}_{d-p},{\bf b}_{p}){\bf 1}_{F_N}({\bf x}_{d-p},{\bf b}_{p}){\bf 1}_{F_N}({\bf x}_{d-p},{\bf y}_{p}){\bf
1}_{F_N}({\bf a}_{d-p},{\bf y}_{p}){\bf 1}_{\{{\bf a}_{d-p}\cap{\bf x}_{d-p}=\emptyset\}}\notag \\ && =\sum_{\gamma=0}^p \gamma!\binom{p}{\gamma}^2
\frac{U_\gamma}{(d!|F_N|)^2},\notag
\end{eqnarray}
where
\begin{eqnarray*}
U_0 &= &{\sum_{\substack{{\bf a}_{d-p},\,{\bf b}_{p} \\{\bf x}_{d-p},\,\,{\bf y}_{p}}}}\!\!{\bf 1}_{F_N}({\bf a}_{d-p},{\bf b}_{p}){\bf 1}_{F_N}({\bf x}_{d-p},{\bf
b}_{p}){\bf 1}_{F_N}({\bf x}_{d-p},{\bf y}_{p}){\bf 1}_{F_N}({\bf a}_{d-p},{\bf y}_{p}){\bf 1}_{\{{\bf a}_{d-p}\cap{\bf x}_{d-p}=\emptyset\}}{\bf 1}_{\{{\bf
b}_{p}\cap{\bf y}_{p}=\emptyset\}} \\ &\leq &|F_N^\#|,
\end{eqnarray*}
and
\begin{eqnarray*}
U_p &= &{\sum_{\substack{{\bf a}_{d-p},\,\, {\bf b}_{p} \\{\bf x}_{d-p}}}}\!\!{\bf 1}_{F_N}({\bf a}_{d-p},{\bf b}_{p}){\bf 1}_{F_N}({\bf x}_{d-p},{\bf b}_{p}){\bf
1}_{\{{\bf a}_{d-p}\cap{\bf x}_{d-p}=\emptyset\}}\\ & \leq & \sum_{{\bf b}_p}\left[ \sum_{{\bf a}_{d-p}} {\bf 1}_{F_N}({\bf a}_{d-p},{\bf b}_{p}) \right]^2\leq
|F_N|\times \max_{j\leq N} |F^*_{N,j}|,
\end{eqnarray*}
and finally, for $\gamma =1,...,p-1$,
\begin{eqnarray*}
U_\gamma &= &{\sum_{\substack{{\bf a}_{d-p},\,{\bf u}_{\gamma},\, {\bf b}_{p-\gamma} \\{\bf x}_{d-p},\,\,{\bf y}_{p-\gamma}}}}\!\!{\bf 1}_{F_N}({\bf a}_{d-p},{\bf
u}_{\gamma}, {\bf b}_{p-\gamma}){\bf 1}_{F_N}({\bf x}_{d-p},{\bf u}_{\gamma}, {\bf b}_{p-\gamma})\times \\ && \quad\quad\quad\quad\quad\quad \times{\bf 1}_{F_N}({\bf
x}_{d-p},{\bf u}_{\gamma}, {\bf y}_{p-\gamma}){\bf 1}_{F_N}({\bf a}_{d-p},{\bf u}_{\gamma}, {\bf y}_{p-\gamma}){\bf 1}_{\{{\bf a}_{d-p}\cap{\bf x}_{d-p}=\emptyset\}}{\bf
1}_{\{{\bf b}_{p-\gamma}\cap{\bf y}_{p-\gamma}=\emptyset\}}\\ &\leq & {\sum_{\substack{{\bf a}_{d-p},\,{\bf u}_{\gamma},\, {\bf b}_{p-\gamma} \\{\bf x}_{d-p},\,\,{\bf
y}_{p-\gamma}}}}\!\!{\bf 1}_{F_N}({\bf a}_{d-p},{\bf u}_{\gamma}, {\bf b}_{p-\gamma}){\bf 1}_{F_N}({\bf x}_{d-p},{\bf u}_{\gamma}, {\bf b}_{p-\gamma})\times \\ &&
\quad\quad\quad\quad\quad\quad \times{\bf 1}_{F_N}({\bf x}_{d-p},{\bf u}_{\gamma}, {\bf y}_{p-\gamma}){\bf 1}_{F_N}({\bf a}_{d-p},{\bf u}_{\gamma}, {\bf
y}_{p-\gamma}){\bf 1}_{\{{\bf b}_{p-\gamma}\cap{\bf y}_{p-\gamma}=\emptyset\}}\\ &\leq & |F_N|^2 \times \|f_N \star_{d-p+1}^{d-p} f_N
\|^2_{\ell^2(\mathbb{N})^{\otimes 2p-1}}.
\end{eqnarray*}
The proof is concluded by using the estimates (\ref{FRLA1}) (for $p=1$) and (\ref{FRLA2}) (for $p=2,...,d-1$). \qed

\medskip

We shall now state a generalization of Theorem \ref{T : BerryEsseenBJ}, providing an explicit upper bound for the normal approximation of multiple integrals defined over
{\sl infinite} sets. This result shows once again that our approach allows to deal directly with the Gaussian approximation of random variables that are functions of the
whole \textsl{infinite} sequence $X$. Consider a sequence of real numbers $\beta = \{\beta_i : i\geq 1\}\in \ell^2(\mathbb{N})$, and define the finite measure on
$\mathbb{N}$
\begin{equation}
m_\beta (A) = \sum_{i\in A} \beta_i^2, \quad A\subset\mathbb{N}.
\end{equation}
We denote by $m_\beta^d$ ($d\geq 2$) the canonical $d$-product measure associated with $m_\beta$. For every $d\geq 2$ and for every set $F\subset \mathbb{N}^d$, we define
the sets $F^\# \subset F\times F$ and $F^*_{j}$, $j=1,2,...$ as before. In particular, $F^*_j$ is the collection of all $(i_1,...,i_d)\in F$ such that $j=i_k$ for some
$k$. For $\beta$ as before, and
 $F\in\Delta_d$ (possibly infinite) such that $m_\beta(F)>0$, we are interested in the normal approximation of the random variable
$$
J(\beta,F) =\frac 1{[d!m_\beta(F)]^{1/2}} \sum_{(i_1,...,i_d)\in F} \beta_{i_1}\cdot\cdot\cdot\beta_{i_d}X_{i_1}\cdot\cdot\cdot X_{i_d}.
$$
The following statement, whose proof (omitted) follows along the lines of that of Theorem 5.4, provides an upper bound for the normal approximation of $J(\beta,F)$.
\begin{prop}\label{P:genJB} Let $Z\sim\mathscr{N}(0,1)$. Under the above notation and assumptions, for every
$h  \in
{\mathcal C}_b^2$
one has
that there exist positive constants $K_1$ and $K_2$, depending uniquely on $d$, $\|h\|_\infty$ and $\|h''\|_\infty$, such that\begin{eqnarray}
&& \big|E[h(Z)] - E[h(J(\beta,F))]\big|\leq K_1 \frac{m_\beta^{2d}(F^\#)^{\frac12}}{m_\beta^{d}(F)} + K_2 \left[\sup_{j\geq
1}\frac{m^d_\beta(F^*_{j})}{m^d_\beta(F)}\right]^{\frac14}.
\end{eqnarray}
\end{prop}
\subsection{Fractional Cartesian products}\label{SS FCProducts}
In this section, we describe an explicit application of Theorem \ref{T : BerryEsseenBJ}. The framework and notation are basically the same as those of Example 1.2 in
\cite{BleiJanson}. Fix integers $d\geq 3$ and $2\leq m\leq d-1$, and consider a collection $\{S_1,...,S_d\}$ of distinct non-empty subsets of $[d] = \{1,...,d\}$ such
that: (i) $S_i \neq \emptyset$, (ii) $\bigcup_i S_i = [d]$, (iii) $|S_i|=m$ for every $i$, (iv) each index $j\in[d]$ appears in exactly $m$ of the sets $S_i$, and (v) the cover $\{S_1,...,S_d\}$ is connected (i.e., it cannot be partitioned into two disjoint partial covers). For every $i=1,...,d$ and every ${\bf y}_d = (y_1,...,y_d)\in \mathbb{N}^d$, we use
the notation $\pi_{S_{i}} {\bf y} = (y_j : j\in S_i)$. Note that the operator $\pi_{S_{i}}$ transforms vectors of $\mathbb{N}^d$ into vectors of $\mathbb{N}^m$. For every
$N\geq d^m$, write $n = \lfloor N^{1/m} \rfloor$, that is, $n$ is the largest integer such that $n\leq N^{1/m}$. Now select a one-to-one map $\varphi$ from $[n]^m$ into
$[N]$, and define
\begin{equation*}
F^*_N = \{(\varphi(\pi_{S_{1}}{\bf k}_d),...,\varphi(\pi_{S_{d}}{\bf k}_d)) : {\bf k}_d\in [n]^d \}\subset [N]^d,
\end{equation*}
(note that, in general, $F^*_N$ is {\sl not} a symmetric set), $F^{**}_N = F^*_N \cap \Delta_N^d,$ and also
\begin{equation}\label{FNfcp}
F_N = {\bf sym}(F^{**}_N),
\end{equation}
where ${\bf sym}(F_N^{**})$ indicates the collections of all vectors ${\bf y}_d=(y_1,...,y_d)\in\mathbb{N}^d$ such that
\begin{equation*}
(y_{\sigma(1)},...,y_{\sigma(d)})\in F^{**}_N
\end{equation*}
for some permutation $\sigma$.
\begin{prop}\label{P : fracCprod}
Let $Z\sim\mathscr{N}(0,1)$, and let  $F_N$ and $\widetilde{S}_N$, $N\geq d^m$, be respectively defined according to (\ref{FNfcp}) and (\ref{MF}). Then,
for every
$h  \in
{\mathcal C}_b^2$,
there exists a constant $K>0$, independent of $N$, such that
$$
\big|E[h(\widetilde{S}_N)] - E[h(Z)]\big|\leq \frac{K}{N^{1/(2m)}}.
$$
\end{prop}
{\it Proof}. The computations contained in \cite[p. 16]{BleiJanson} can be straightforwardly adapted to our setting, so that we deduce that the sequence $\{F_N : N\geq
d^m\}$ has \textsl{combinatorial dimension} $\alpha = d/m$. Recall that this means that there exist finite constants $0<Q_2 <Q_1 < \infty$ (independent of $N$) such that:
(a) $|F_N|\geq Q_2 N^\alpha$, and (b) for every $A_1,...,A_d \subset [N]$, $|F_N\cap (A_1\times\cdot\cdot\cdot A_d)|\leq Q_1 (\max_{1\leq j\leq d} |A_j|)^{\alpha}$.
Thanks to Theorem \ref{T : BerryEsseenBJ}, to conclude the proof it is therefore sufficient to check that, as $N\rightarrow\infty$,
\begin{equation}\label{proofsteps}
\max_{j\leq N}|F^*_{N,j}| = O(N^{\alpha-1}) \,\,\,\,\,\, \text{and}\,\,\,\,\,\, |F^\#_{N}| = O(N^{2\alpha-1/m}).
\end{equation}
Start by observing that every element $(z_1,...,z_d)$ of $F_N$ has the form
\begin{equation}\label{zed}
(z_1,...,z_d) = (\varphi(\pi_{S_{\sigma(1)}}{\bf k}_d),...,\varphi(\pi_{S_{\sigma(d)}}{\bf k}_d)),
\end{equation}
where ${\bf k}_d = (k_1,...,k_d)\in[n]^d$ and $\sigma$ is a permutation of $[d]$. Since $\varphi$ is one-to-one, it follows that, for every $j\leq N$, there are at most
$d!d \times n^{d-m}$ elements of the set $F_{N,j}$. To see this, just observe that, every $(z_1,...,z_d)\in F_{N,j}$ is completely specified by the following three
elements: (i) a permutation $\sigma$ of $d$, (ii) the index $a =1,...,d$ such that $z_a = j$, and (iii) the values of those coordinates $k_b$ such that $b\not\in
S_{\sigma(a)}$. Since $n\leq N^{1/m}$ by construction, one immediately obtains the first relation in (\ref{proofsteps}).

\noindent To prove the second part of (\ref{proofsteps}), we shall first show that, if $((z_1,...,z_d),(z'_1,...,z'_d))\in F_N^{\#}$ are such
 that $(z_1,...,z_d)$ is as in (\ref{zed}), and
 $(z'_1,...,z'_d) = (\varphi(\pi_{S_{\sigma'(1)}}{\bf k'}_d),...,\varphi(\pi_{S_{\sigma'(d)}}{\bf k'}_d))$ (${\bf k'}_d\in [n]^d$) then
 ${\bf k}_d$ and ${\bf k}'_d$ must have $m$ coordinates in common. The definition of $F_N^{\#}$ implies indeed that,
for such a vector $((z_1,...,z_d),(z'_1,...,z'_d))$, there exists
$$
(u_1,...,u_d) = (\varphi(\pi_{S_{\rho(1)}}{\bf i}_d),...,\varphi(\pi_{S_{\rho(d)}}{\bf i}_d)) \in F_N,
$$
such that, for some $p=1,...,d-1$: (a) there exist indices $a_1,...,a_p$ and $b_1,...,b_p$ such that
 $\pi_{S_{\rho(a_i)}}{\bf i}_d = \pi_{S_{\sigma(b_i)}}{\bf k}_d$ for every $i=1,...,p$, and (b) there exist
indices $v_1,...,v_{d-p}$ and $w_{1},...,w_{d-p}$ such that $\pi_{S_{\rho(v_i)}}{\bf i}_d = \pi_{S_{\sigma'(w_i)}}{\bf k'}_d$ for every $i=1,...,d-p$. Note that then we
necessarily have $\{a_1,...,a_p,v_1,...,v_{d-p}\} = \{b_1,...,b_p,w_1,...,w_{d-p}\} = [d]$. By connectedness, there
 exists at least one $q\in [d]$ such that $i_q$ (i.e., the $q$th coordinate of ${\bf i}_d$) belongs both to one of
the sets $\pi_{S_{\rho(a_i)}}{\bf i}_d$ and to one of the sets $\pi_{S_{\rho(v_i)}}{\bf i}_d$. 
This last property implies that there
 exists a constant $L$, independent of $N$, such that $|F^\#_N|\leq L n^{2d-1}\leq L N^{2\alpha -1/m}$. This concludes the proof.
\qed

\begin{rem}{\rm
Note that the combinatorial dimension $\alpha= d/m$, as appearing in the proof of Proposition \ref{P : fracCprod}, is an index of interdependence between the coordinates
of the sets $F_N$. See Ch. XIII in \cite{BleiBook} for more details on this point. }
\end{rem}

\subsection{Beyond the Rademacher case: a question by Blei and Janson}\label{SS : questionBJ}
Now we go back to the framework and notation of Section \ref{SSGenBJ}, so that, in particular, the sequence $\widetilde{S}_N$, $N\geq 1$, is defined according to
(\ref{MF}).
For every $N$ define $\widetilde{S}^G_N$ to be the random variable obtained from (\ref{MF}) by replacing the sequence $X$ with a i.i.d. Gaussian sequence
$G=\{G_i : i\geq 1\}$, where each $G_i$ has mean zero and
unit variance. A natural question, which has been left open by Blei and Janson in \cite[Remark
4.6]{BleiJanson}, is
 whether under the conditions (\ref{BJneg1})--(\ref{BJneg2}) the sequence $\widetilde{S}^G_N$, $N\geq 1$, converges in law towards a standard
Gaussian distribution. Note that this problem could be tackled by a direct computation, based for instance on \cite{NouPec_ptrf} or \cite{NP}. However, the results of
this paper, combined with those of \cite{MOO}, allow to elegantly deduce a more general result, which we provide in the forthcoming statement. In what follows, we write $V
= \{V_i : i\geq 1\}$ to indicate a centered i.i.d. sequence, with
unit variance and such that $E|V_1|^3= \eta <\infty$ (note that the results of \cite{MOO} would allow
to obtain similar results in even more general frameworks, but we do not look for generality here). We also denote by $\widetilde{S}^V_N$ the random variable obtained
from (\ref{MF}) by replacing $X$ with $V$.

\begin{prop} Under the above notation and assumptions, consider a
function
$h  \in
{\mathcal C}_b^3$.
Then, there exist two universal constants $B_1$ and $B_2$, depending uniquely on $d$, $\eta$,
$\|h\|_\infty$, $\|h''\|_\infty$ and $\|h'''\|_\infty$, such that, for every $N\geq 1$ and for $Z\sim\mathscr{N}(0,1)$,
\begin{equation}\label{BJboundgauss}
 \big|E[h(Z)] - E[h(\widetilde{S}^V_N)]\big|
 \leq B_1 \frac{|F^\#_{N}|^{\frac12}}{|F_{N}|}
+B_2 \left[\max_{j\leq N}\frac{|F^*_{Nj}|}{|F_{N}|}\right]^{\frac14}.
\end{equation}
In particular, if (\ref{BJneg1})--(\ref{BJneg2}) take place, then $\widetilde{S}^V_N$ converges in law towards $Z$.
\end{prop}
{\it Proof.} One has that
$$
\big|E[h(Z)] - E[h(\widetilde{S}^V_N)]\big|\leq \big|E[h(Z)] - E[h(\widetilde{S}_N)]\big| +\big|E[h(\widetilde{S}_N)] - E[h(\widetilde{S}^V_N)]\big|,
$$
and the conclusion is obtained by combining Theorem \ref{T : BerryEsseenBJ} with the fact that, according to \cite[Theorem 3.18, case H2]{MOO}, there exists a constant
$Q$, depending only on $\|h'''\|_\infty$ and $\eta$, such that
$$
\big|E[h(\widetilde{S}_N)] - E[h(\widetilde{S}^V_N)]\big|\leq Q \left[\max_{j\leq N}\frac{|F^*_{Nj}|}{|F_{N}|}\right]^{\frac12}.
$$
\qed

\vskip1cm

\noindent {\bf Acknowledgments.} This paper originates from all three authors attending the ``Workshop on Stein's method''
at the Institute of Mathematical Sciences of the
National University of Singapore,
from March 31 until April 4, 2008. We warmly thank A. Barbour, L. Chen, K.-P.
Choi and A. Xia for their kind hospitality and generous support. We are grateful to J. Dedecker for useful discussions.

\section*{Appendix: Some technical proofs}\label{Appendix}

\subsection*{Proof of Lemma \ref{L : estimates}} 
In what follows, for $k\geq 1$, we shall write ${\bf a}_k$ to indicate the generic vector ${\bf a}_k=(a_1,...,a_k) \in \N^k$. We start by proving the first point. To do
this, just write, using the Cauchy-Schwarz inequality,
\begin{eqnarray*}
\|f\star_r^lg\|^2_{\ell^2(\N)^{\otimes(n+m-r-l)}}&=& \sum_{{\bf i}_{m+n-r-l}}       f\star_r^lg\,({\bf i}_{m+n-r-l})^2\\ &=& \sum_{{\bf i}_{n-r}} \sum_{{\bf j}_{m-r}}
\sum_{{\bf k}_{r-l}}  \left( \sum_{{\bf a}_l} f({\bf i}_{n-r},{\bf k}_{r-l},{\bf a}_l) g({\bf j}_{m-r},{\bf k}_{r-l},{\bf a}_l) \right)^2\\ &\leq& \sum_{{\bf i}_{n-r}}
\sum_{{\bf k}_{r-l}} \sum_{{\bf a}_l}  f^2({\bf i}_{n-r},{\bf k}_{r-l},{\bf a}_l) \sum_{{\bf j}_{m-r}}\sum_{{\bf b}_l}  g^2({\bf j}_{m-r},{\bf k}_{r-l},{\bf b}_l)\\
&\leq& \sum_{{\bf i}_{n-r}} \sum_{{\bf k}_{r-l}} \sum_{{\bf a}_l}  f^2({\bf i}_{n-r},{\bf k}_{r-l},{\bf a}_l) \sum_{{\bf j}_{m-r}}\sum_{{\bf l}_{r-l}} \sum_{{\bf b}_l}
g^2({\bf j}_{m-r},{\bf l}_{r-l},{\bf b}_l)\\ &=&\|f\|^2_{\ell^2(\N)^{\otimes n}} \|g\|^2_{\ell^2(\N)^{\otimes m}}.
\end{eqnarray*}
The first part of the second point is obtained by writing
\begin{eqnarray*}
\max_j \left(\sum_{{\bf b}_{n-1}} f^2(j,{\bf b}_{n-1})\right)^2 &\leq& \sum_j \left(\sum_{{\bf b}_{n-1}} f^2(j,{\bf b}_{n-1})\right)^2\\ &=& \|f\star_{n}^{n-1}
f\|^2_{\ell^2(\N)}\\ &\leq& \max_{j} \sum_{{\bf b}_{n-1}} f^2(j,{\bf b}_{n-1}) \times \sum_{j'} \sum_{{\bf b}'_{n-1}} f^2(j,{\bf b}_{n-1})\\ &=& \max_{j} \sum_{{\bf
b}_{n-1}} f^2(j,{\bf b}_{n-1}) \times \|f\|^2_{\ell^2(\N)^{\otimes n}}.
\end{eqnarray*}

For the second part, we have, by applying (in order) the Cauchy-Schwarz inequality, as well as the previous estimate,
\begin{eqnarray*}
&&\|f\star_l^{l-1}g\|^2_{\ell^2(\N)^{\otimes(n+m-2l+1)}}\\ &=&\sum_j\sum_{{\bf k}_{n-l}}\sum_{{\bf l}_{m-l}} \left(\sum_{{\bf i}_{l-1}} f(j,{\bf i}_{l-1},{\bf
k}_{n-l})g(j,{\bf i}_{l-1},{\bf l}_{m-l})\right)^2=\sum_j \|f(j,\cdot)\star_{l-1}^{l-1}g(j,\cdot)\|^2_{\ell^2(\N)^{\otimes(n+m-2l)}}\\ &\leq&\sum_j\sum_{{\bf k}_{n-l}}
\sum_{{\bf i}_{l-1}} f^2(j,{\bf i}_{l-1},{\bf k}_{n-l}) \sum_{{\bf l}_{m-l}} \sum_{{\bf j}_{l-1}} g^2(j,{\bf j}_{l-1},{\bf l}_{m-l})\\ &\leq& \max_j \sum_{{\bf b}_{n-1}}
f^2(j,{\bf b}_{n-1}) \times \sum_{{\bf i}_{m}} g^2({\bf i}_m) = \max_j \sum_{{\bf b}_{n-1}} f^2(j,{\bf b}_{n-1}) \times \|g\|^2_{\ell^2(\N)^{\otimes m}}\\ &=&
\sqrt{\max_j \left(\sum_{{\bf b}_{n-1}} f^2(j,{\bf b}_{n-1})\right)^2} \times \|g\|^2_{\ell^2(\N)^{\otimes n}}\\ &\leq& \sqrt{\sum_j \left(\sum_{{\bf b}_{n-1}} f^2(j,{\bf
b}_{n-1})\right)^2} \times \|g\|^2_{\ell^2(\N)^{\otimes m}}=\|f\star_{n-1}^{n-1} f\|_{\ell^{2}(\N)} \times \|g\|^2_{\ell^2(\N)^{\otimes m}}.
\end{eqnarray*}

Finally, for the third part, just observe that
$$
\|f\star_1^0 f\|^2_{\ell^2(\N)^{\otimes (2n-1)}} =\sum_j \sum_{{\bf i}_{n-1}}\sum_{{\bf j}_{n-1}} f(j,{\bf i}_{n-1})f(j,{\bf j}_{n-1}) =\sum_j \left(\sum_{{\bf i}_{n-1}}
f(j,{\bf i}_{n-1})\right)^2=\|f\star_{n}^{n-1} f\|_{\ell^2(\N)}^2
$$
and, for $2\leq l\leq n$:
\begin{eqnarray*}
&&\|f\star_{l}^{l-1}f\|^2_{\ell^2(\N)^{\otimes(2n-2l+1)}}\\ &=&\sum_{i=j}\sum_{{\bf k}_{n-l}}\sum_{{\bf l}_{n-l}} \left(\sum_{{\bf i}_{l-1}} f(i,{\bf i}_{l-1},{\bf
k}_{n-l})f(j,{\bf i}_{l-1},{\bf l}_{n-l})\right)^2\\ &\leq&\sum_{({\bf a}_{n-l+1},{\bf b}_{n-l+1})\in\Delta_{2n-2l+2}^c} \left(\sum_{{\bf i}_{l-1}} f({\bf i}_{l-1},{\bf
a}_{n-l+1})f({\bf i}_{l-1},{\bf b}_{n-l+1})\right)^2\\ &=&\|f\star_{l-1}^{l-1} f\times {\bf 1}_{\Delta^c_{2n-2l+2}}\|^2_{\ell^2(\N)^{\otimes (2n-2l+2)}}.
\end{eqnarray*}

\subsection*{Proof of Lemma \ref{L : convLemma}} 

We have, by bilinearity and using the first point of Lemma \ref{L : estimates},
\begin{eqnarray*}
&&\|f_k\star_r^r g_k - f\star_r^r g\|_{\ell^2(\N)^{\otimes(m+n-2r)}} =\|f_k\star_r^r (g_k-g) + (f_k-f)\star_r^r g\|_{\ell^2(\N)^{\otimes(m+n-2r)}}\\ && \leq\|f_k\star_r^r
(g_k-g)\|_{\ell^2(\N)^{\otimes(m+n-2r)}} + \|(f_k-f)\star_r^r g\|_{\ell^2(\N)^{\otimes(m+n-2r)}}\\ &&\leq\|f_k\|_{\ell^2(\N)^{\otimes n}} \|g_k-g\|_{\ell^2(\N)^{\otimes
m}} +\|f_k-f\|_{\ell^2(\N)^{\otimes n}} \|g\|_{\ell^2(\N)^{\otimes m}} \underset{k\to\infty}{\longrightarrow} 0.
\end{eqnarray*}

\subsection*{Proof of Proposition \ref{Product formula}}\label{A : ProofProduct}
We assume, without loss of generality, that $m\leq n$. We first prove (\ref {s}) for functions $f$ and $g$ with support, respectively, in $\Delta _{n}^{N}$ and $\Delta _{m}^{N}$, for some finite $N$. One has that
\begin{equation*}
J_{n}\left( f\right) J_{m}\left( g\right) =\int_{\Delta _{n}^{N}\times \Delta _{m}^{N}}\left\{ f\star _{0}^{0}g\right\} d\mu _{X}^{\otimes n+m}=\int_{\Delta
_{n}^{N}\times \Delta _{m}^{N}}\left\{ f\otimes g\right\} d\mu _{X}^{\otimes n+m}.
\end{equation*}
For $r=0,...,m$, we denote by $\Pi _{r} \left( n,m\right) $ the set of all partitions $\pi$ of $\left\{ 1,....,n+m\right\} $ composed of

\begin{description}
\item[i)] exactly $r\ $blocks of the type $\left\{ i_{1},i_{2}\right\} $, with $1\leq i_{1}\leq n$ and $n+1\leq i_{2}\leq n+m$,

\item[ii)] exactly $n+m-2r$ singletons.
\end{description}

For instance: an element of $\Pi _{2} \left( 3,2\right) $ is the partition
\begin{equation}
\pi =\left\{ \left\{ 1,4\right\} ,\left\{ 2,5\right\} ,\left\{ 3\right\} \right\} ;  \label{sss}
\end{equation}
the only element of $\Pi _{0} \left( n,m\right) $ is\ the partition $\pi =\left\{ \left\{ 1\right\} ,\left\{ 2\right\} ,...,\left\{ n+m\right\} \right\} $ composed of all
singletons; an element of $\Pi _{3} \left( 4,4\right) $ is\ the partition
\begin{equation}
\pi =\left\{ \left\{ 1,5\right\} ,\left\{ 2,6\right\} ,\left\{ 3,7\right\} ,\left\{ 4\right\} ,\left\{ 8\right\} \right\} . \label{ssss}
\end{equation}
It is easily seen that $\Pi _{r} \left( n,m\right) $ contains exactly $r!\binom{n}{r}\binom{m}{r}$ elements (to specify an element of $\Pi _{r} \left( n,m\right) $, first select $r$ elements of $\left\{ 1,...,n\right\} $, then select $r$ elements in $\left\{ n+1,...,n+m\right\} $ , then build a bijection between the two selected $r$-sets). For every $ r=0,...,m$ and every $\pi \in \Pi _{r} \left( n,m\right) $, we write $B_{n,m}^{N}\left( \pi \right) $ to denote the subset of $\left\{ 1,...,N\right\} ^{n+m}$ given by
\begin{equation*}
\left\{ \left( i_{1},...,i_{n},i_{n+1},...,i_{n+m}\right) \in \left\{ 1,...,N\right\} ^{n+m}:i_{j}=i_{k}\text{ iff }j\text{ and }k\text{ are in the same block of }\pi \right\}
\end{equation*}
(note the ``if and only if'' in the definition). For instance, for $\pi$ as in (\ref{sss}), an element of $B_{3,2}^{3}\left( \pi \right) $ is $\left( 1,2,3,1,2\right) $;
for $\pi$ as in (\ref{ssss}), an element of $B_{4,4}^{5}\left( \pi \right) $ is $\left( 1,2,3,4,1,2,3,5\right) $. The following two facts can be easily checked;

\begin{description}
\item[A)]
\begin{equation*}
\Delta _{n}^{N}\times \Delta _{m}^{N}=\bigcup\limits_{r=0}^{m}\bigcup\limits_{\pi \in \Pi _{r} \left( n,m\right) }B_{n,m}^{N}\left( \pi \right) ,
\end{equation*}
where the unions are disjoint;

\item[B)] For every $r=0,...,m$, and every $\pi \in \Pi _{r} \left( n,m\right) $,
\begin{eqnarray*}
\int_{B_{n,m}^{N}\left( \pi \right) }\left\{ f\star _{0}^{0}g\right\} d\mu _{X}^{\otimes n+m} &=&\int_{\Delta _{n+m-2r}}\left\{ f\star _{r}^{r}g\right\} d\mu
_{X}^{\otimes n+m-2r}.
\end{eqnarray*}
(note that the last expression does not depend on the partition $\pi$, but only on the class $\Pi _{r} \left( n,m\right) $).
\end{description}
\noindent It follows that
\begin{eqnarray*}
J_{n}\left( f\right) J_{m}\left( g\right) &=&\int_{\Delta _{n}^{N}\times \Delta _{m}^{N}}\left\{ f\star _{0}^{0}g\right\} d\mu _{X}^{\otimes n+m} \\ &=&\sum_{r=0}^{m}\sum\limits_{\pi ^{\left( r\right) }\in \Pi _{r}^{\ast }\left( n,m\right) }\int_{B_{n,m}^{N}\left( \pi ^{\left( r\right) }\right) }\left\{ f\star _{0}^{0}g\right\} d\mu _{X}^{\otimes n+m} \\ &=&\sum_{r=0}^{m}r!\binom{n}{r}\binom{m}{r} \int_{\Delta _{n+m-2r}}\left\{ f\star _{r}^{r}g\right\} d\mu _{X}^{\otimes n+m} \\ &=^{(*)}& \sum_{r=0}^{m}r!\binom{n}{r}\binom{m}{r} \int_{\Delta _{n+m-2r}}\left\{ \widetilde{f\star _{r}^{r}g}\right\} d\mu _{X}^{\otimes n+m} \\ &=&\sum_{r=0}^{m}r!\binom{n}{r}\binom{m}{r}J_{n+m-2r}\left[ \left( \widetilde{f\star _{r}^{r}g}\right) \mathbf{1}_{\Delta _{n+m-2r}}\right] ,
\end{eqnarray*}%
which is the desired conclusion. We stress that the equality $(*)$ has been obtained by using the symmetry of the measure $\mu _{X}^{\otimes n+m}$. The result for general
$f,g$ is deduced by an approximation argument and Lemma \ref{L : convLemma}. This concludes the proof. \qed
\end{document}